\documentclass[11pt,a4paper]{article}
\usepackage[utf8]{inputenc}
\usepackage[T1]{fontenc}
\usepackage{geometry}
\usepackage{amsmath, amssymb, amsthm}
\usepackage{graphicx}
\usepackage[labelformat=parens,font=footnotesize,justification=centering]{subcaption}
\usepackage{hyperref}
\usepackage{xcolor}
\usepackage{mathrsfs}
\usepackage{enumitem}
\usepackage{bm}
\usepackage{mathtools}
\usepackage{comment}
\newtheorem{theorem}{Theorem}[section]
\newtheorem{lemma}[theorem]{Lemma}
\newtheorem{proposition}[theorem]{Proposition}
\newtheorem{conjecture}[theorem]{Conjecture}
\newtheorem{corollary}[theorem]{Corollary}

\newtheorem{remark}[theorem]{Remark}
\newtheorem{definition}[theorem]{Definition}

\newcommand{\tx}{\tilde{x}}
\newcommand{\ta}{\tilde{a}}
\newcommand{\tw}{\tilde{w}}

\newcommand{\bydef}{\,\stackrel{\mbox{\tiny\textnormal{\raisebox{0ex}[0ex][0ex]{def}}}}{=}\,} 
\newcommand{\norm}[1]{\left\|#1\right\|}
\newcommand{\abs}[1]{|#1|}
\newcommand{\nuL}{\ell^1_\nu}

\newcommand{\im}{\mathrm{i}}
\newcommand{\C}{\mathbb{C}}

\newcommand{\N}{\mathbb{N}}
\newcommand{\R}{\mathbb{R}}

\newcommand{\Z}{\mathbb{Z}}

\newcommand{\cp}{\bm{\ast}}

\newcommand{\ba}{\bar{a}}

\newcommand{\bx}{\bar{x}}

\newcommand{\akt}[1]{{#1}}

\graphicspath{{figs}}

\begin{document}

\title{Star-shaped periodic solutions of a system of delay differential equations: A computer-assisted proof}
\date{\today}
\author{%
Shoya~Motonaga \thanks{Department of Artificial Intelligence, Kyushu Institute of Technology, 680-4 Kawazu, Iizuka, Fukuoka 820-8502, Japan (\texttt{motonaga.shoya640@mail.kyutech.jp})
}%
\and
Akitoshi~Takayasu \thanks{Institute of Systems and Information Engineering, University of Tsukuba, 1-1-1 Tennodai, Tsukuba, Ibaraki 305-8573, Japan (\texttt{takitoshi@risk.tsukuba.ac.jp})
}
}
\maketitle
%
%\tableofcontents
% 

\begin{abstract}
In this paper, we consider a planar system of nonlinear delay differential equations and prove the
existence of the star-shaped periodic solutions conjectured by Hara (Nonlinear Analysis: Theory,
Methods \& Applications, 30:4555--4563 (1997)) for seven parameter values. These solutions cannot be
obtained by a direct application of the standard Fourier spectral method, since the
corresponding linearized operator is degenerate. From numerical observations of the Fourier
coefficients, we derive an ansatz that writes a star-shaped periodic solution as a periodic
profile in a rotating frame. Restricting the profile to real Fourier coefficients, we show that the
resonance between the delay and the period of the profile fixes the rotation frequency exactly.
We then validate a zero of the resulting real system with the radii polynomial approach and
interval arithmetic, with bounds that cover all Fourier modes. For each parameter value, we also
determine the fundamental period and the winding number of the solution.
\end{abstract}
{\raggedright\noindent\textbf{Keywords:} delay differential equations, star-shaped periodic solutions,
computer-assisted proof, rigorous numerics, Fourier spectral method, radii polynomial approach\par}


% ****************************************************************************************
% Section 1 Introduction 
% ****************************************************************************************

\section{Introduction}

Delay differential equations (DDEs) describe systems whose rate of change depends not only on
the current state but also on the past states.  They arise in a wide variety of applications,
including population dynamics, neural networks, optics, engineering systems and feedback control.
The initial value problem for a DDE is posed with an initial history, and its solution defines a
dynamical system whose phase space is an infinite-dimensional space of history functions; see,
for example, \cite{hale1993introduction, diekmann1995delay, wu1996theory}.  As a consequence, the
mathematical analysis of DDEs is often substantially more difficult than that of ordinary
differential equations (ODEs), and the analysis of nonlinear DDEs remains challenging even for
periodic solutions.  Periodic and quasi-periodic oscillations are nevertheless among the most
fundamental behaviors observed in applications, and establishing their existence is an important
problem both mathematically and practically.

The infinite-dimensional phase space also changes the geometry of solutions.  For an autonomous
ODE, the uniqueness of solutions to the initial value problem implies that a trajectory in the
phase space never intersects itself.  For a DDE, the same holds for the orbit of the history in
the infinite-dimensional space of history functions.  However, the projection of this orbit onto
a finite-dimensional subspace may intersect itself, and so may the solution curve itself in the
Euclidean space.  In other words, even a planar DDE can have periodic solutions whose curves
cross themselves, which is impossible for a planar autonomous ODE.

In the present paper, we consider the following system of DDEs:
\begin{align}\label{eq:DDE}
	\dot{\bm{x}}(t)=-\alpha\left(1-|\bm{x}(t)|^{2}\right) R(\theta) \bm{x}(t-1),
\end{align}
where $\alpha\in\mathbb{R}$,
\[
\bm{x}= \begin{bmatrix}
	x\\y
\end{bmatrix},~ R(\theta)=\begin{bmatrix}
	\cos\theta & -\sin\theta\\\sin\theta & \cos\theta
\end{bmatrix}~\left(|\theta|<\frac{\pi}{2}\right),~ |\bm{x}|^2=x^2+y^2.
\]
Equivalently, introducing the complex coordinate $z(t)=x(t)+\im y(t)$,
\begin{align}\label{eq:DDE_complex}
\dot{z}(t)=-\alpha \left(1-|z(t)|^2\right) e^{\im\theta} z(t-1).
\end{align}
Hara and his collaborators \cite{P1,P2,P3} have introduced the system \eqref{eq:DDE} by
incorporating the rotation matrix $R(\theta)$ in the plane into the scalar delay differential
equations
\[
\dot{x}(t)=-\alpha{\left(1+x(t)\right)} x(t-1),\quad\dot{x}(t)=-\alpha{\left(1-x^{2}(t)\right)} x(t-1),
\]
They have proved part of the behavior of its solutions and have observed further behavior
numerically.
First, Hara and Sugie~\cite{P1} have shown that the zero solution is asymptotically stable for sufficiently small $\alpha$.
\begin{theorem}[Hara--Sugie\cite{P1}]\label{thm:P1}
	The zero solution of $\dot{\bm{y}}(t) = -\alpha R(\theta) \boldsymbol{y}(t-1)$ is asymptotically stable if and only if $0<\alpha<{\pi}/{2}-|\theta|$.
\end{theorem}
Next, Hara~\cite{P2} has shown that, as $\alpha$ increases, a Hopf bifurcation occurs and a nontrivial periodic solution emerges.
\begin{theorem}[Hara\ \cite{P2}]\label{thm:P2}
	If $\alpha>{\pi}/{2}-|\theta|>0$, then \eqref{eq:DDE} has non-constant periodic solutions  such that
	\begin{align*}
		\begin{cases}
			x(t) = r\cos\omega t\\
			y(t) = r\sin\omega t
		\end{cases}~(\theta >0),\quad
		\begin{cases}
			x(t) = r\cos(-\omega t)\\
			y(t) = r\sin(-\omega t)
		\end{cases}~(\theta <0),
	\end{align*}
	where
	\[
	r=\sqrt{1-\frac{1}{\alpha}\left(\frac{\pi}{2}-|\theta|\right)}, \quad \omega=-\left(\frac{\pi}{2}-|\theta|\right).
	\]
\end{theorem}
Motivated by whether \eqref{eq:DDE} possesses periodic solutions other than the circular orbit $(CO)$ of Theorem~\ref{thm:P2}, Hara~\cite{P2} has also carried out mathematical considerations, which have led to the following conjecture:
\begin{conjecture}[Hara\,\cite{P2}]\label{conj:hara}
	There exists $\alpha_{0}>{\pi}/{2}-|\theta|>0$ such that
	\begin{enumerate}
		\item if $\alpha_{0} \geq \alpha>{\pi}/{2}-|\theta|$, then $(CO)$ is the only nonconstant {stable} periodic solution of \eqref{eq:DDE}.
		\item if $\alpha>\alpha_{0}$, then $(CO)$ is an {unstable} periodic solution. Moreover, when $\theta={\pi}/{n}$ and the initial function $(\phi(t), \psi(t))$ satisfies $0<\phi^{2}(0)+\psi^{2}(0)<1$:
		\begin{itemize}
			\item If $n$ is a {rational number}, the solution orbit of \eqref{eq:DDE} asymptotically converges to a {star-shaped (or polygonal) periodic solution} with $p$-vertices.
			\item If $n$ is an irrational number, each solution orbit densely fills out the annular region centered at the origin with radii $\sin \theta$ and 1.
		\end{itemize}
	\end{enumerate}
\end{conjecture}

\begin{figure}[htbp]
\centering
\begin{subfigure}[t]{0.24\textwidth}\centering
\includegraphics[width=\linewidth]{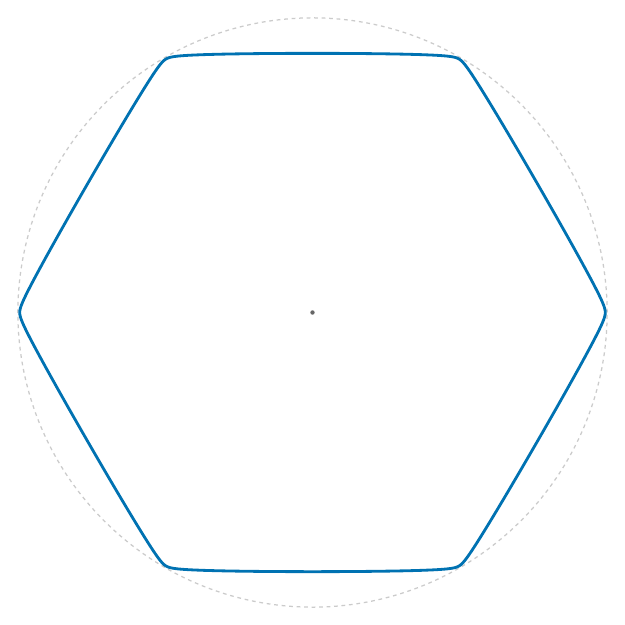}
\caption{$\theta=\pi/3$, $\alpha=6$\\ $(p,m)=(6,1)$}\label{fig:star_pi3}
\end{subfigure}\hfill
\begin{subfigure}[t]{0.24\textwidth}\centering
\includegraphics[width=\linewidth]{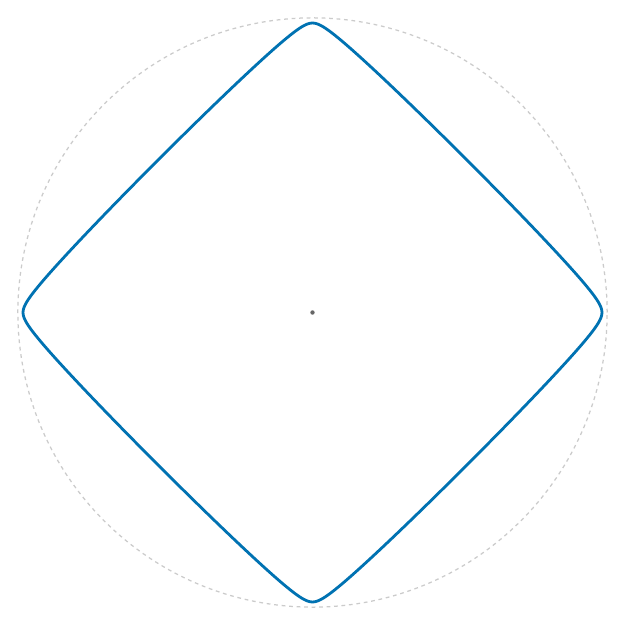}
\caption{$\theta=\pi/4$, $\alpha=4$\\ $(p,m)=(4,1)$}\label{fig:star_pi4}
\end{subfigure}\hfill
\begin{subfigure}[t]{0.24\textwidth}\centering
\includegraphics[width=\linewidth]{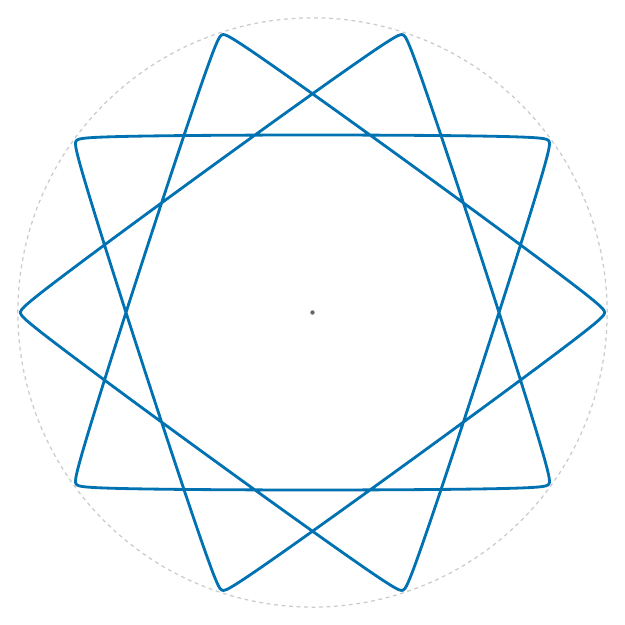}
\caption{$\theta=\pi/5$, $\alpha=4$\\ $(p,m)=(10,3)$}\label{fig:star_pi5}
\end{subfigure}\hfill
\begin{subfigure}[t]{0.24\textwidth}\centering
\includegraphics[width=\linewidth]{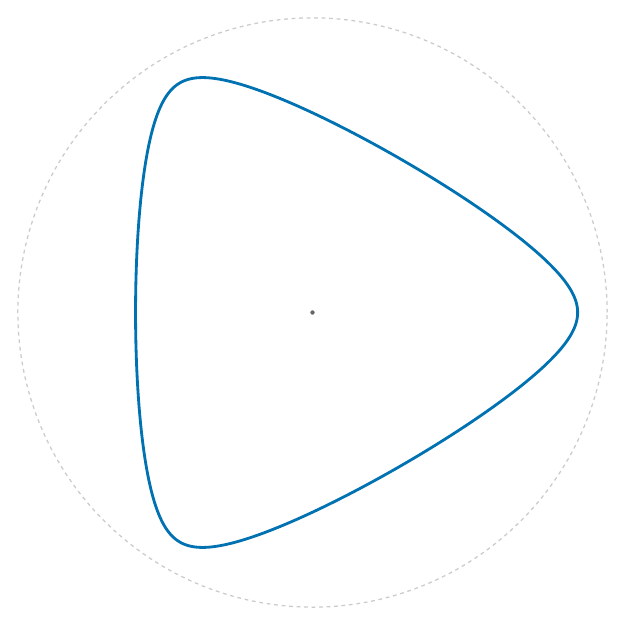}
\caption{$\theta=\pi/6$, $\alpha=13/5$\\ $(p,m)=(3,1)$}\label{fig:star_pi6}
\end{subfigure}

\vspace{3mm}
\begin{subfigure}[t]{0.24\textwidth}\centering
\includegraphics[width=\linewidth]{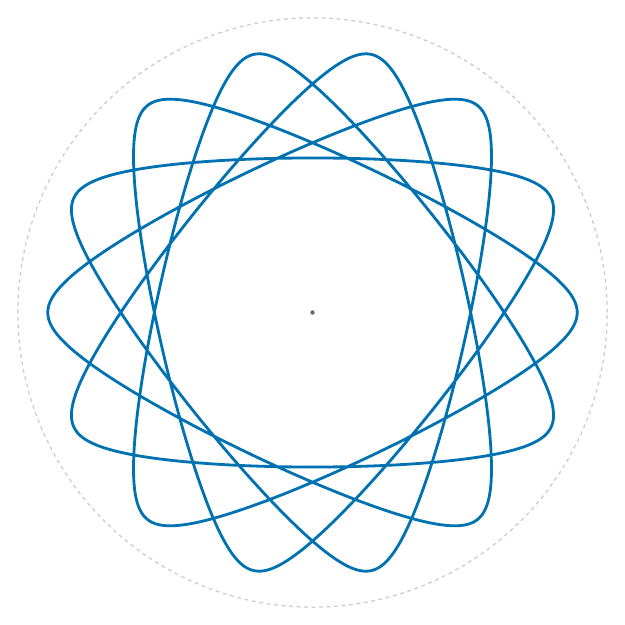}
\caption{$\theta=\pi/7$, $\alpha=5/2$\\ $(p,m)=(14,5)$}\label{fig:star_pi7}
\end{subfigure}\hspace{0.012\textwidth}
\begin{subfigure}[t]{0.24\textwidth}\centering
\includegraphics[width=\linewidth]{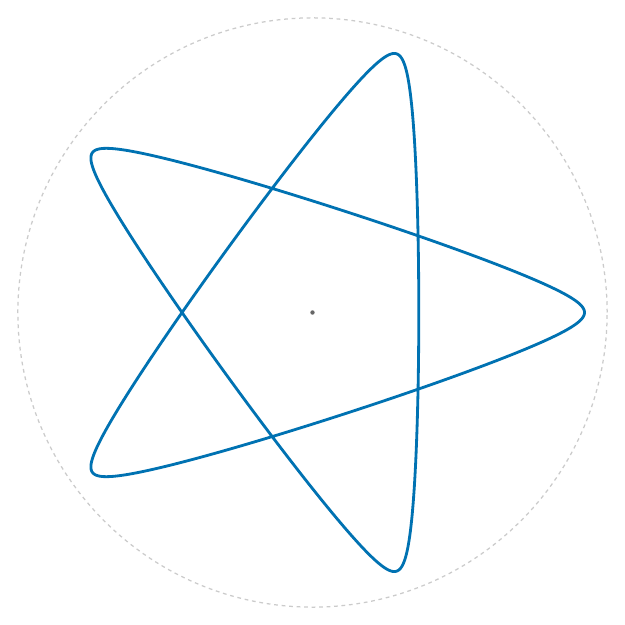}
\caption{$\theta=\pi/10$, $\alpha=5/2$\\ $(p,m)=(5,2)$}\label{fig:star_pi10}
\end{subfigure}\hspace{0.012\textwidth}
\begin{subfigure}[t]{0.24\textwidth}\centering
\includegraphics[width=\linewidth]{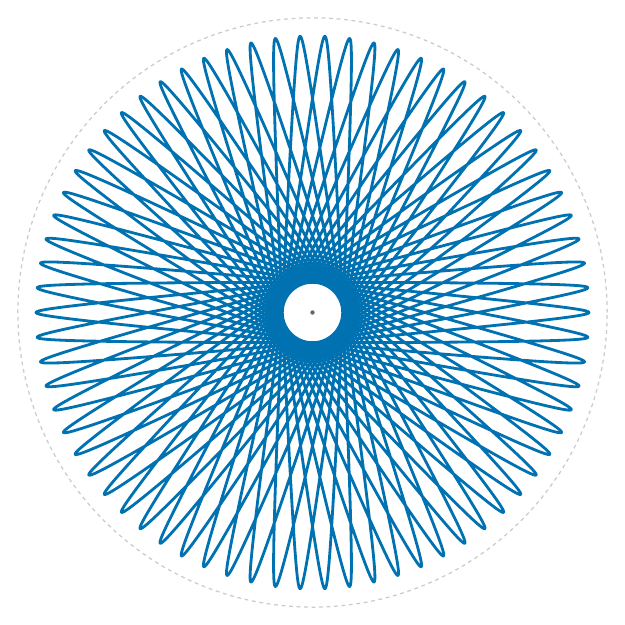}
\caption{$\theta=\pi/35$, $\alpha=5/2$\\ $(p,m)=(70,33)$}\label{fig:star_pi35}
\end{subfigure}
\caption{Star-shaped periodic solutions to \eqref{eq:DDE} for the seven parameter pairs of
Theorem~\ref{thm:main}. Each panel shows the curve $t\mapsto z(t)$ in the complex plane over one
fundamental period $L=2p$. The curves are computed from the numerical approximations used in the
computer-assisted proof, and the validated error bounds of Section~\ref{sec:cap} are far below the
resolution of the plots. The solution rotates $m$ times around the origin and is invariant under
rotation by $2\pi/p$, as stated in Theorem~\ref{thm:main}. All panels use the same scale, and the
dashed curve is the unit circle.}
\label{fig:star_orbits}
\end{figure}

In particular, Hara~\cite{P2} has numerically observed periodic solutions whose curves draw
star-shaped figures in the plane when $\theta/\pi$ is rational (see Figure~\ref{fig:star_orbits}).  These results were based on
numerical observations, and no mathematical proof of the existence of such star-shaped periodic
solutions has been given since the conjecture was formulated.  In addition to the nonlinearity,
the infinite-dimensional nature of \eqref{eq:DDE} makes the problem difficult.  This is our main
motivation for this study: we prove the existence of star-shaped periodic solutions by \emph{computer-assisted proofs}.

% In connection with this conjecture, Hara and Sakata~\cite{P3} proved the existence of star-shaped periodic solutions for the
% related equation with piecewise constant arguments
% \[
% \dot{\bm{x}}(t)=-\alpha\left(1-|\bm{x}(t)|^2\right)R(\theta)\bm{x}(\lfloor t\rfloor),
% \]
% where $\lfloor\cdot\rfloor$ denotes the floor function. This equation precisely differs from the
% delay system \eqref{eq:DDE}, so those existence results do not settle the latter
% problem.
% Conjecture~\ref{conj:hara} itself remains unresolved.
% % was studied in \cite{P3} , and the same solution behavior as in Conjecture~\ref{conj:hara} was proved due to the role of the floor function.

% \begin{theorem}[Hara--Sakata~\cite{P3}]\label{thm:P3}
% There exists $\alpha_{0}>0$ such that\vspace*{-2mm}
% \begin{enumerate}
% 	\item if $0<\alpha \leqslant \alpha_{0}$, then there does not exist any nonconstant periodic solution.\vspace*{-2mm}
% 	\item if $\alpha>\alpha_{0}$ and $\theta / \pi\in\mathbb{Q}$, then there exist star-shaped periodic solutions.\vspace*{-2mm}
% 	\item if $\alpha>\alpha_{0}$ and $\theta / \pi\not\in\mathbb{Q}$, then each solution orbit with $\left\|\bm{x}\left(0\right)\right\|<1$ densely fills out an annular region centered at the origin.
% \end{enumerate}
% \end{theorem}

In recent years, computer-assisted proofs have become a powerful framework for studying nonlinear
differential equations.  They have been developed for ODEs, DDEs, and partial differential equations; see, for example, the monograph
\cite{nakao2019numerical} and the surveys \cite{vandenberg2015rigorous, gomezserrano2019computer}.
For DDEs, computer-assisted proofs of periodic solutions have been
obtained by several groups.  A representative example is Wright's
equation, for which Fourier-based Newton--Kantorovich methods have led to the proofs of the
conjectures of Wright and Jones \cite{lessard2010recent, vandenberg2018proof, jaquette2019proof}.  These methods have been
extended to systems with several delays \cite{kiss2012computational, vandenberg2022general}.
For the Mackey--Glass equation, another representative example of DDEs, the existence of periodic orbits has also been proved with rigorous integrators
\cite{szczelina2018algorithm, szczelina2024high} and with a combination of analysis and computer
assistance \cite{bartha2021stable}.

Among these approaches, Newton-type methods based on Fourier spectral expansions are one of the
standard techniques for proving the existence of periodic solutions.  One reformulates the problem
as a zero-finding problem in a Banach space of Fourier coefficients, computes an approximate zero numerically, and then applies the Radii Polynomial Theorem, whose hypotheses are verified with
interval arithmetic.  The method combines high numerical accuracy with a rigorous a posteriori
validation, and it proves existence and local uniqueness of a periodic solution near the
approximate one.

The argument requires that the linearized operator at the approximate solution has a bounded
inverse.  When this operator is degenerate, the validation fails.  Such situations frequently
occur in the presence of symmetries, resonances or small divisors associated with periodic or
quasi-periodic dynamics. Such difficulties have been overcome in specific settings, for
example in computer-assisted proofs of invariant tori with quasi-periodic dynamics
\cite{celletti2007kam, figueras2017rigorous, capinski2020computer}, but there is currently no
general-purpose framework for handling them.
The system \eqref{eq:DDE} is a case in point.  A previous numerical study by the second author
\cite{Takayasu-RIMS} computed approximate star-shaped periodic solutions with a Fourier spectral
method and observed a characteristic structure: their periods appear to be of the form $2p$
with $p\in\mathbb{N}$, and only the Fourier coefficients with indices $pl\pm m$ ($l\in\mathbb{Z}$)
appear to be numerically nonzero for some $m\in\mathbb{N}$.  The same study also suggested an obstruction to a computer-assisted proof.  Star-shaped periodic solutions seem to appear only
when $\theta/\pi$ is rational, which suggests a singular dependence on the parameter, and when the
Fourier spectral method is applied to \eqref{eq:DDE} directly, the linearized operator of the
zero-finding problem was observed to become numerically non-invertible.

The main purpose of this paper is to remove this degeneracy.  Revisiting the numerical
observations of \cite{Takayasu-RIMS}, we find that a star-shaped periodic solution decomposes into
a rotation with constant speed and a $2$-periodic deformation.  This leads to the following
ansatz, which reflects the resonance structure of the solutions:
\begin{quote}
\textbf{Ansatz.}
The corresponding (quasi-)periodic solutions of \eqref{eq:DDE_complex} in Conjecture~\ref{conj:hara} are of the form
\begin{align}\label{eq:main-anzatz}
z(t)=w(t)e^{\im\mu t}=r(t)e^{\im \rho(t)}e^{\im \mu t},
\end{align}
where $w(t)$, $r(t)$, $\rho(t)$ are some 2-periodic functions and $\mu\in\R$ is a constant parameter.
\end{quote}
The motivation for this ansatz comes from the numerical observations above and from the
theoretical analysis in Section~\ref{sec:ansatz}.  The ansatz solution has three periodic factors,
$r(t)$ ($2$-periodic), $e^{\im\rho(t)}$ ($2$-periodic) and $e^{\im\mu t}$ ($2\pi/\mu$-periodic).  A solution of this form is quasi-periodic in general, and it becomes periodic if $\mu$
is a rational multiple of $\pi$; the converse holds under a nondegeneracy condition (see
Proposition~\ref{prop:dich_star}).  In other words, the ansatz separates the rational rotation,
which causes the singular dependence on $\theta$, from the $2$-periodic profile $w$.  Substituting
the ansatz into \eqref{eq:DDE_complex} gives a complex-valued DDE for $w$ with the parameter $\mu$,
and the problem reduces to finding a parameter $\mu$ and a $2$-periodic solution $w$ of this equation.  The reduced problem no longer has the degeneracy described above. We determine
$\mu$ algebraically by a symmetry reduction and then validate $w$ with the Radii Polynomial
Theorem.

We now return to our target, the star-shaped periodic solutions of \eqref{eq:DDE}, and make
precise which solutions we prove to exist.
The noncircular orbits observed by Hara~\cite{P2}
motivate a class of solutions characterized by a time translation followed by a rotation. For the unit-delay system
\eqref{eq:DDE}, we adopt the following definition:

\begin{definition}[Star-shaped periodic solution]\label{def:star_shaped}
Let $p,m$ be coprime positive integers with $0<m<p$. A solution
$\bm{x}(t)=(\operatorname{Re}z(t),\operatorname{Im}z(t))$ of \eqref{eq:DDE}, defined for
all $t\in\R$, is a \emph{star-shaped periodic solution of type $(p,m)$} if $|z(t)|$ is
nonconstant and
\begin{equation}\label{eq:star_shaped_symmetry}
z(t+2)=e^{-2\pi\im m/p}z(t)\qquad(t\in\R).
\end{equation}
\end{definition}

Iterating \eqref{eq:star_shaped_symmetry} gives $z(t+2p)=z(t)$, so $2p$ is a period.
Since $p$ and $m$ are coprime, the orbit is invariant under rotation by $2\pi/p$.
The nonconstant modulus excludes the circular solutions of Theorem~\ref{thm:P2}.
We use the term \emph{star-shaped} in this spatio-temporal sense throughout the paper;
the definition does not prescribe a number of radial extrema or a polygonal shape.
The integer $m$ specifies the rotation in \eqref{eq:star_shaped_symmetry}.
Neither the fundamental period nor the winding number is part of the definition.
We establish both for the solutions in the following theorem.

\vspace*{1mm}
The main contribution of the present paper is to give computer-assisted proofs of the
existence of star-shaped periodic solutions to \eqref{eq:DDE} for seven parameter values, shown
in Figure~\ref{fig:star_orbits}. This
establishes the existence of the star-shaped periodic solutions observed numerically in
\cite{P2}, which had remained unproved. Our main theorem is the following:

\begin{theorem}[Existence of star-shaped periodic solutions]\label{thm:main}
Let $\theta=\pi/n$ and let $(n,\alpha)$ be one of the seven pairs
\[
(3,6),\quad (4,4),\quad (5,4),\quad (6,13/5),\quad (7,5/2),\quad (10,5/2),\quad (35,5/2).
\]
Then \eqref{eq:DDE} has a star-shaped periodic solution of type $(p,m)$ in the sense of
Definition~\ref{def:star_shaped}, where $m/p=(n-2)/(2n)$ in lowest terms. The corresponding
pairs are
\[
(p,m)=(6,1),\ (4,1),\ (10,3),\ (3,1),\ (14,5),\ (5,2),\ (70,33),
\]
respectively. The solution can be written as $z(t)=w(t)e^{\im\mu t}$, where $w$ is a
$2$-periodic complex-valued function and $\mu=\theta-\pi/2=-m\pi/p$.
Its fundamental period is $L=2p$. It avoids the origin and has signed winding number
$-m$, that is, $m$ clockwise turns over one fundamental period.
\end{theorem}

Theorem~\ref{thm:main} gives an existence result at specific rational values of
$\theta/\pi$ for the class specified by Definition~\ref{def:star_shaped}. This class
captures the rotational structure motivating Conjecture~\ref{conj:hara}.
The proof combines a symmetry reduction that fixes the rotation frequency exactly with
interval estimates that validate a zero of the resulting Fourier system. Two nonzero
Fourier coefficients imply that the modulus is nonconstant and determine the fundamental
period. A further interval bound gives the winding number.
The stability and the attraction of these solutions and the instability threshold $\alpha_0$
for the circular orbit $(CO)$ remain open.

% Theorem~\ref{thm:main} gives a partial answer to part~2 of Conjecture~\ref{conj:hara}: it
% settles the existence of the star-shaped periodic solutions on the rational branch
% $\theta/\pi\in\mathbb{Q}$ at specific parameter values.  The remaining assertions of the
% conjecture, that is, the threshold $\alpha_0$, the stability statements and the claim that
% generic solutions asymptotically converge to these orbits, are of a dynamical nature and are not
% addressed in this paper.  Section~\ref{subsec:dich_ledger} records what is and what is not
% proved.  The proof of Theorem~\ref{thm:main} is obtained in Sections~\ref{sec:dde}
% and~\ref{sec:dde} by a careful blend of a Fourier spectral discretization, a
% Newton--Kantorovich type theorem with a frozen preconditioner, interval arithmetic and a
% symmetry analysis of the resonance between the delay and the period.  The sufficient conditions
% of the Newton--Kantorovich type theorem are rigorously checked by numerical computations based
% on interval arithmetic.

The rest of the paper is organized as follows.
{
In Section~\ref{sec:ansatz}, we explain how Definition~\ref{def:star_shaped} and the ansatz
are derived from numerical observations of the Fourier coefficients of the star-shaped periodic
solutions, and we give the DDE in a rotating frame under the ansatz.
% In Section~\ref{sec:dde}, we seek profiles of Fourier representation for the star-shaped periodic solutions, and reduce the existence of the periodic solutions to that of zeros of an induced map on a suitable Fourier space.
% In Section~\ref{sec:cap}, we give the  {computer-assisted} proof of the  {star-shaped periodic} solutions of our system. The theoretical details of the associated tail multiplier, bounded map, and estimates for the corresponding bounds as well as their numerical values are also provided.
In Section~\ref{sec:dde}, we write the DDE in the rotating frame as a zero-finding problem on a
weighted space of Fourier coefficients. Restricting the profile to real Fourier coefficients, we
show that the rotation frequency is determined exactly as $\mu=\theta-\pi/2$, and we obtain a real
system for the profile. We also show that a zero of this system gives a star-shaped periodic
solution and determine its fundamental period.
In Section~\ref{sec:cap}, we validate a zero of the real system with the Radii Polynomial Theorem,
derive the bounds $Y$, $Z_0$ and $Z_1$ \akt{of this theorem for the zero-finding map,} including all Fourier modes, and prove
Theorem~\ref{thm:main}. We close the paper with concluding remarks on the scope of the result
and future work.
}

\section{Derivation of the ansatz and the DDE in a rotating frame}\label{sec:ansatz}
In this section, we explain how Definition~\ref{def:star_shaped} is derived from numerical
observations of the Fourier coefficients of the solutions.  We formulate
the observations as conditions (i) and (ii) below, derive a first form of the solution from
condition~(i) in Lemma~\ref{lem:(i)}, refine it with condition~(ii) in Lemma~\ref{lem:(ii)}, and
arrive at the symmetry in Definition~\ref{def:star_shaped}.  The role of
Sections~\ref{subsec:ansatz_i} and~\ref{subsec:ansatz_ii} is mainly to explain where the definition and
the ansatz come from; the computer-assisted proof in Sections~\ref{sec:dde} and~\ref{sec:cap} is
built on the formulation given in Section~\ref{sec:numerical_periodic}.

We start from the following numerical observations.
\begin{enumerate}[label=(\roman*),leftmargin=2.5em]
\item The period $L$ of the star-shaped periodic solutions of \eqref{eq:DDE} is of the form $2p\ (p\in\mathbb{N})$~\cite{Takayasu-RIMS}.
\item Write the Fourier expansions of the periodic solutions of \eqref{eq:DDE} as
\[
x(t)=\sum_{k\in \Z} {\xi_k} e^{\im \omega k t},\qquad y(t)=\sum_{k\in \Z} {\eta_k} e^{\im \omega k t},\qquad \xi_k,\eta_k\in\mathbb{C}
\]
with $\omega=2\pi/L=\pi/p$. In the numerical solutions computed in this study, for some
$m\in\mathbb{N}$, the complex Fourier coefficients $\xi_k+\im\eta_k$ of $z=x+\im y$ appear to
vanish unless $k\equiv-m\pmod p$. In particular, $\xi_k$ and $\eta_k$ vanish unless
$k=pl\pm m\ (l\in\mathbb{Z})$, as observed in~\cite{Takayasu-RIMS}.
\end{enumerate}
In other words, the Fourier coefficients of these periodic solutions are supported on a
restricted set of indices.
These conditions lead to the ansatz \eqref{eq:main-anzatz} and to
Definition~\ref{def:star_shaped}, as we show in the following subsections.

{Throughout the paper, $c^\ast$ denotes the complex conjugate of $c\in\C$ and we reserve the bar for approximate solutions.}

\subsection{Implication of condition (i): polar form in the complex plane}\label{subsec:ansatz_i}
We first use only condition~(i).  Writing a solution of \eqref{eq:DDE_complex} in polar
coordinates, we obtain a weaker form of the ansatz, in which $r$ and $\rho$ are $2p$-periodic
rather than $2$-periodic.
\begin{lemma}\label{lem:(i)}
Assume that
\[
z(t)=r(t)e^{\im \lambda(t)}
\]
is a periodic solution of \eqref{eq:DDE_complex} with period $2p$ $(p\in\mathbb{N})$, where $r(t)>0$ and $\lambda:\mathbb{R}\to\mathbb{R}$ are real-valued smooth functions.
Then $r$ is $2p$-periodic and
\[
\lambda(t)=\rho(t)+\pi\frac{m}{p} t
\]
for some $m\in\mathbb{Z}$ and a $2p$-periodic function $\rho$.
\end{lemma}
\begin{proof}
	Since $z$ is $2p$-periodic, so is $r=|z|$.
    {Using $z(t)=r(t)e^{\im \lambda(t)}$, we have
    \[
    \dot{z}(t)=\dot{r}(t)e^{\im \lambda(t)}+\im\dot{\lambda}(t)r(t)e^{\im \lambda(t)}.
    \]
    Dividing by $z(t)$, we see that
    \[
    \frac{\dot z}{z}=\frac{\dot r}{r}+\im\dot\lambda,
    \]
    so that $\dot\lambda(t)=\operatorname{Im}\bigl(\dot z(t)/z(t)\bigr)$.
    By the $2p$-periodicity of $z$, we deduce that $\dot\lambda(t)$ is $2p$-periodic.
	Note that $\lambda(t)$ itself need not be $2p$-periodic.}
	
	Put $\mu:=\frac{1}{2p}\int_0^{2p} \dot{\lambda}(s) ds$ and $\rho(t)=\int_0^t (\dot{\lambda}(s)-\mu)ds+\lambda(0)$.
	Then
	\begin{align*}
	\lambda(t)&=\lambda(0)+\int_0^t \dot{\lambda}(s) ds\\
	&=\lambda(0)+\mu t+\int_0^t (\dot{\lambda}(s)-\mu) ds\\
	&=\mu t+\rho(t).
	\end{align*}
	By the $2p$-periodicity of $\dot{\lambda}$, we see that
	\[
	\rho(t+2p)-\rho(t)=\int_t^{t+2p} (\dot{\lambda}(s)-\mu) ds=\int_0^{2p} (\dot{\lambda}(s)-\mu) ds=0,
	\]
	i.e., $\rho$ is $2p$-periodic.
	Since $z,~r$, and $\rho$ are $2p$-periodic, $e^{\im \mu t}$ is also $2p$-periodic,
	which implies that $2p\mu=2\pi m$ for some $m\in\mathbb{Z}$.
	Therefore, we obtain $\mu=\pi \frac{m}{p}$.
\end{proof}

\subsection{Implication of condition (ii): analysis in the $(x,y)$-plane}\label{subsec:ansatz_ii}
{Next}, we show that the vanishing condition (ii) for the Fourier coefficients imposes
a strong constraint on the form of the solution{, which leads to the ansatz}.
\begin{lemma}\label{lem:(ii)}
 Assume that conditions (i) and (ii) hold, and
 let $z=x+\im y$ be a $2p$-periodic solution of \eqref{eq:DDE_complex} with $z(t)\neq0$ for all
 $t\in\R$.
  Then $z$ is of the form
 \begin{align*}
z(t)=r(t)e^{\im\rho(t)}e^{\im\pi\frac{m'}{p}t}
\end{align*}
with 2-periodic functions $r,\rho$ and an integer $m'$ with $m'\equiv-m\pmod p$.
Moreover, $r$, $\rho$ and $m'$ coincide with those in Lemma~\ref{lem:(i)}, up to adding a
constant in $2\pi\Z$ to $\rho$, and $z(t+2)=e^{-2\pi\im m/p}z(t)$ for all $t\in\R$.
\end{lemma}
It is worth noting that the factor $e^{\im\pi\frac{m'}{p}t}$ is $2p$-periodic and rotates
$|m'|$ times around the unit circle during $t\in[0,2p]$, while $r e^{\im\rho}$ deforms the orbit
with period~$2$; the resulting orbits are shown in Figure~\ref{fig:star_orbits}.
\begin{proof}
By condition~(ii), we have
\begin{align*}
x(t)&=\sum_{k\in\mathbb{Z}} \xi_k e^{\im\frac{\pi}{p}kt}=\left(\sum_{l\in\mathbb{Z}} \xi_{pl+m} e^{\im\pi lt}\right) e^{\im\pi\frac{m}{p}t}
+\left(\sum_{l'\in\mathbb{Z}} \xi_{pl'-m} e^{\im\pi l't}\right) e^{-\im\pi\frac{m}{p}t}.
\end{align*}

{Since $x(t)$ is real-valued, {we have $\xi_{k}=\xi^\ast_{-k}$ for $k\in\mathbb{Z}$.}
%SH:", where $c^\ast$ denotes the complex conjugate of $c\in\C$. We use this notation throughout the paper, and reserve the bar for approximate solutions." moved to the opening of section 2.
Thus we have
\[
\sum_{l'\in\mathbb{Z}} \xi_{pl'-m} e^{\im\pi l't}=\sum_{l'\in\mathbb{Z}} \xi_{-pl'-m} e^{-\im\pi l' t}=\sum_{l'\in\mathbb{Z}} \xi^\ast_{pl'+m} e^{-\im\pi l' t}.
\]
Then we obtain
\begin{align*}
x(t)&=\left(\sum_{l\in\mathbb{Z}} \xi_{pl+m} e^{\im\pi lt}\right) e^{\im\pi\frac{m}{p}t}
+\left(\sum_{l'\in\mathbb{Z}}  {\xi^\ast_{pl'+m}} e^{-\im\pi l't}\right) e^{-\im\pi\frac{m}{p}t}\\
&=\sum_{l\in\mathbb{Z}} (\xi_{pl+m}+ {\xi^\ast_{pl+m}})\cos\pi\left(l+\frac{m}{p}\right)t+
\im\sum_{l\in\mathbb{Z}} (\xi_{pl+m}- {\xi^\ast_{pl+m}})\sin\pi\left(l+\frac{m}{p}\right)t\\
&=2\sum_{l\in\mathbb{Z}} \mathrm{Re}(\xi_{pl+m})\cos\pi\left(l+\frac{m}{p}\right)t
-2\sum_{l\in\mathbb{Z}} \mathrm{Im}(\xi_{pl+m})\sin\pi\left(l+\frac{m}{p}\right)t\\
&=2\left(\sum_{l\in\mathbb{Z}} \left(\mathrm{Re}(\xi_{pl+m})\cos\pi lt- \mathrm{Im}(\xi_{pl+m})\sin\pi lt\right)\right)\cos\pi\frac{m}{p}t\\
&\qquad -2\left(\sum_{l\in\mathbb{Z}} \left(\mathrm{Re}(\xi_{pl+m})\sin\pi lt+\mathrm{Im}(\xi_{pl+m})\cos\pi lt\right)\right)\sin\pi\frac{m}{p}t.
\end{align*}

Similarly, we obtain
\begin{align*}
y(t)&=2\left(\sum_{l\in\mathbb{Z}} \left(\mathrm{Re}(\eta_{pl+m})\cos\pi lt- \mathrm{Im}(\eta_{pl+m})\sin\pi lt\right)\right)\cos\pi\frac{m}{p}t\\
&\qquad -2\left(\sum_{l\in\mathbb{Z}} \left(\mathrm{Re}(\eta_{pl+m})\sin\pi lt+\mathrm{Im}(\eta_{pl+m})\cos\pi lt\right)\right)\sin\pi\frac{m}{p}t.
\end{align*}}%
Therefore, there exist 2-periodic functions $\varphi_x, \psi_x,\varphi_y,\psi_y$
such that
\begin{align*}
\begin{pmatrix}
x(t)\\
y(t)
\end{pmatrix}
=\begin{pmatrix}
\varphi_x(t)\cos(\pi\frac{m}{p}t)+\psi_x(t)\sin(\pi\frac{m}{p}t)\\
\varphi_y(t)\cos(\pi\frac{m}{p}t)+\psi_y(t)\sin(\pi\frac{m}{p}t)
\end{pmatrix}
=\Phi(t)\begin{pmatrix}
\cos(\pi\frac{m}{p}t)\\
\sin(\pi\frac{m}{p}t)
\end{pmatrix}
\end{align*}
with the 2-periodic matrix
\begin{align*}
\Phi(t)=\begin{pmatrix}
\varphi_x(t)& \psi_x(t)\\
\varphi_y(t)& \psi_y(t)
\end{pmatrix}.
\end{align*}

In terms of $A_x(t)\bydef\sum_{l\in\Z}\xi_{pl+m}e^{\im\pi lt}$ and
$A_y(t)\bydef\sum_{l\in\Z}\eta_{pl+m}e^{\im\pi lt}$, the above expansions read
$\varphi_x=2\operatorname{Re}A_x$, $\psi_x=-2\operatorname{Im}A_x$,
$\varphi_y=2\operatorname{Re}A_y$ and $\psi_y=-2\operatorname{Im}A_y$.
By condition~(ii), $\xi_k+\im\eta_k=0$, i.e., $\eta_k=\im\xi_k$, for $k\equiv m\pmod p$, so that
$A_y=\im A_x$, and hence
$\varphi_y=\psi_x$ and $\psi_y=-\varphi_x$, that is,
\[
\Phi(t)=\begin{pmatrix}
\varphi_x(t)& \varphi_y(t)\\
\varphi_y(t)& -\varphi_x(t)
\end{pmatrix}.
\]
Therefore,
\[
z(t)=x(t)+\im y(t)
=\bigl(\varphi_x(t)+\im\varphi_y(t)\bigr)\Bigl(\cos\bigl(\pi\tfrac{m}{p}t\bigr)-\im\sin\bigl(\pi\tfrac{m}{p}t\bigr)\Bigr)
=\bigl(\varphi_x(t)+\im\varphi_y(t)\bigr)e^{-\im\pi\frac{m}{p}t}.
\]
In particular, $z(t+2)=e^{-2\pi\im m/p}z(t)$, since $\varphi_x$ and $\varphi_y$ are $2$-periodic.
Since $z(t)\neq0$, we can write $\varphi_x(t)+\im\varphi_y(t)=r(t)e^{\im\rho_0(t)}$ with $r(t)>0$
and a continuous function $\rho_0$.  As $\varphi_x$ and $\varphi_y$ are $2$-periodic, so are $r$
and $e^{\im\rho_0}$.  Hence $\rho_0(t+2)-\rho_0(t)\in2\pi\Z$, and by continuity
$\rho_0(t+2)-\rho_0(t)=2\pi j$ for an integer $j$ independent of $t$.  Put
\[
\rho(t)\bydef\rho_0(t)-\pi jt,\qquad m'\bydef pj-m.
\]
Then $\rho(t+2)-\rho(t)=2\pi j-2\pi j=0$, so $\rho$ is $2$-periodic, and
$z(t)=r(t)e^{\im\rho(t)}e^{\im\pi jt}e^{-\im\pi\frac{m}{p}t}=r(t)e^{\im\rho(t)}e^{\im\pi\frac{m'}{p}t}$
with $m'\equiv-m\pmod p$.

It remains to compare with Lemma~\ref{lem:(i)}.  The function
$\lambda(t)\bydef\rho(t)+\pi\frac{m'}{p}t$ is a continuous argument of $z$, and $\rho$ is
$2p$-periodic, so this is a representation of the form in Lemma~\ref{lem:(i)}.  Such a
representation is unique in the following sense.  A continuous argument of $z$ is determined up
to an additive constant in $2\pi\Z$, and the integer in Lemma~\ref{lem:(i)} is determined by
$\lambda(2p)-\lambda(0)=2\pi m'$, since $\rho$ is $2p$-periodic.  Hence $m'$ and $\rho$ coincide
with those in Lemma~\ref{lem:(i)}, the latter up to a constant in $2\pi\Z$, which completes the
proof.
\end{proof}

Lemma~\ref{lem:(ii)} is the end point of the derivation.  Its conclusion is exactly the
symmetry \eqref{eq:star_shaped_symmetry} in Definition~\ref{def:star_shaped}, and it gives the
ansatz \eqref{eq:main-anzatz} with $2$-periodic $r,\rho$ and $\mu=\pi m'/p$.  The choice
$\mu=-\pi m/p$ used in Section~\ref{sec:numerical_periodic} differs from it by $\pi j$, and the
factor $e^{\im\pi jt}$ is $2$-periodic; both represent the same solution.  Since $m'$ is the
winding number of $z$ over $[0,2p]$ and $j$ is that of $re^{\im\rho_0}$ over one period, the
relation $m'=pj-m$ is the formula in Remark~\ref{rem:dich_table_correction}.  The other two requirements of the definition are not derived from the
observations.  The nonconstant modulus excludes the circular solutions of
Theorem~\ref{thm:P2}, and the conditions $0<m<p$ and $\gcd(p,m)=1$ only normalize the pair
$(p,m)$, since \eqref{eq:star_shaped_symmetry} depends only on $m/p$ modulo $1$.  This is how we
arrived at Definition~\ref{def:star_shaped}.

\subsection{The DDE in a rotating frame}\label{sec:numerical_periodic}
Lemmas~\ref{lem:(i)} and~\ref{lem:(ii)} show how the numerical observations (i) and (ii) lead
to Definition~\ref{def:star_shaped}.  Since these conditions are observations, the derivation
does not prove that such solutions exist, and it determines $\mu$ only up to the form
$\mu\in\pi\Z/p$.  From now on, we take Definition~\ref{def:star_shaped} as the starting point.

Here, using the ansatz \eqref{eq:main-anzatz}, we derive the equation in a rotating frame and express the symmetry in
Definition~\ref{def:star_shaped} as periodicity of the profile. We work with the complex form
\eqref{eq:DDE_complex} of \eqref{eq:DDE} and assume $\alpha>0$ from now on.  We seek solutions of
the form \eqref{eq:main-anzatz}, $z(t)=w(t)e^{\im\mu t}$, where $w$ is a complex-valued
$2$-periodic profile and $\mu\in\R$ is a parameter.  We use $L$ for the fundamental period
of $z$.
Since $|z(t)| = |w(t)|$, $z(t-1) = w(t-1)\,e^{\im\mu(t-1)}$ and
$\dot{z}(t) = \bigl(\dot{w}(t) + \im\mu\, w(t)\bigr) e^{\im\mu t}$, substituting the ansatz
into \eqref{eq:DDE_complex} and dividing both sides by $e^{\im\mu t}$, we obtain the $w$-DDE
\begin{align}\label{eq:w_dde}
  \dot{w}(t) = -\alpha\, e^{\im(\theta - \mu)}\bigl(1 - |w(t)|^2\bigr)\, w(t-1) - \im\mu\, w(t).
\end{align}
We write $\beta \bydef \alpha\, e^{\im(\theta-\mu)}$ and note that $|w|^2 = w {w^\ast}$.

The parameter $\mu$ in \eqref{eq:main-anzatz} is not a parameter of the original DDE \eqref{eq:DDE_complex}.
It is a decomposition parameter that specifies the rotating frame.
Given a solution $z$ of \eqref{eq:DDE_complex}, the function $w(t) = e^{-\im\mu t} z(t)$ solves \eqref{eq:w_dde} for every $\mu \in \R$, and the ansatz \eqref{eq:main-anzatz} asks that $\mu$ be chosen so that $w$ is $2$-periodic.
From the ansatz, the relation between successive periods is
\begin{align}\label{eq:relative_periodicity}
  z(t + 2) = z(t)\, e^{2\im\mu}.
\end{align}
Hence, $z$ is a \emph{relative periodic orbit} with multiplier $e^{2\im\mu}$.
If $\mu/\pi\in\mathbb{Q}$, then $z$ is periodic. The converse requires a
nondegeneracy condition: a one-mode profile can represent a periodic circular orbit even
when $\mu/\pi$ is irrational. Such a condition is given in
Proposition~\ref{prop:dich_star}, which also gives the exact period.

For a solution of type $(p,m)$, choose $\mu=-m\pi/p$. Then
\eqref{eq:star_shaped_symmetry} is equivalent to $2$-periodicity of
$w(t)=e^{-\im\mu t}z(t)$. Conversely, a $2$-periodic solution of \eqref{eq:w_dde} with
this frequency and a nonconstant modulus gives a star-shaped periodic solution by
Definition~\ref{def:star_shaped}.

It remains to construct a $2$-periodic profile with nonconstant modulus. We do this in
the symmetric class $w(-t)= {w(t)^\ast}$. The next section shows that, for
nonequilibrium solutions in this class, we may take $\mu=\theta-\pi/2$ by choosing a
representative of the rotating frame. This exact choice will provide the rational
rotation relation required by Definition~\ref{def:star_shaped}.

\section{Fourier representation and symmetry reduction}\label{sec:dde}

We first express the rotating-frame equation in a weighted Fourier space. We then
restrict the construction to profiles satisfying $w(-t)= {w(t)^\ast}$, for which the
Fourier coefficients are real. On this subspace, the resonance between the delay and
the profile period determines the rotation frequency modulo $\pi$. Choosing a
representative gives a real system whose zeros solve the full DDE.

\subsection{Fourier space and the differential equation}\label{subsec:dde_zero}

We seek profiles of period  {$2$}. For $\nu\ge1$, define the space of Fourier coefficients as
\begin{equation}\label{eq:ell1nu_complex}
\nuL\bydef\left\{a\in\C^{\Z}:\norm a_\nu\bydef
\sum_{k\in\Z}|a_k|\nu^{|k|}<\infty\right\}.
\end{equation}
An element $a\in\nuL$ represents the continuous function
$w(t)=\sum_{k\in\Z}a_ke^{\im {\pi} kt}$. The delay and complex conjugation act by
\begin{equation}\label{eq:shift_dagger}
(Sa)_k\bydef e^{-\im {\pi} k}a_k=(-1)^ka_k,
\qquad (a^\dagger)_k\bydef {a_{-k}^\ast}.
\end{equation}
With the discrete convolution $(a\cp b)_k=\sum_{j\in\N}a_jb_{k-j}$, the cubic term $|w|^2w(t-1)$ has
coefficient sequence
\begin{equation}\label{eq:g_cubic}
g(a)\bydef a\cp a^\dagger\cp Sa.
\end{equation}
The space $\nuL$ is a Banach algebra, with
$\norm{a\cp b}_\nu\le\norm a_\nu\norm b_\nu$, and $S$ and $\dagger$ are isometries.
Thus $\norm{g(a)}_\nu\le\norm a_\nu^3$ holds. We regard this space as a real Banach space
when differentiating, since $\dagger$ is conjugate-linear.

Substituting this series into \eqref{eq:w_dde}, we obtain the mode equations $G_k(a,\mu)=0$,
$k\in\Z$, where
\begin{equation}\label{eq:Gk}
G_k(a,\mu)\bydef\lambda_k a_k-\beta[g(a)]_k,
\qquad \lambda_k\bydef\im {\pi} k+\im\mu+\beta(-1)^k {.}
\end{equation}
The derivative in the equation is unbounded on $\nuL$. Its natural domain is denoted by
\[
\mathcal D_\nu\bydef
\left\{a\in\nuL:\sum_{k\in\Z}|k||a_k|\nu^{|k|}<\infty\right\}.
\]
Nevertheless, any $a\in\nuL$ satisfying every mode equation $G_k(a,\mu)=0$ belongs
to this domain. Indeed,
\[
\im {\pi} k a_k=-\im\mu a_k-\beta(Sa)_k+\beta[g(a)]_k
\]
has an $\ell^1_\nu$ right-hand side. Its Fourier series therefore has an absolutely
summable derivative and defines a classical solution of \eqref{eq:w_dde}.
Conversely, a profile in this domain solves \eqref{eq:w_dde} only if all its mode
equations vanish. Section~\ref{sec:cap} uses a bounded preconditioned map to validate
these equations on the sequence space.

\subsection{The symmetric class and the rotation frequency}\label{subsec:dich_splitting}

A profile satisfies $w(-t)= {w(t)^\ast}$ if and only if its Fourier coefficients
are real. We seek solutions in the closed real subspace
\begin{equation}\label{eq:real_space}
X\bydef\nuL(\R)=\{a\in\nuL:a_k\in\R\text{ for all }k\in\Z\},
\qquad \norm a_X=\norm a_\nu.
\end{equation}
This restriction is part of the existence construction, rather than a requirement
of Definition~\ref{def:star_shaped}. It suffices to find solutions in this class to
prove the existence asserted in Theorem~\ref{thm:main}.

Since the profile period $2$ is twice the delay $1$, the resonance identity is
\begin{equation}\label{eq:dich_resonance}
e^{-\im\pi k}=(-1)^k\in\R,
\qquad S^2=I.
\end{equation}
Thus the delay shift preserves real coefficient sequences. The following splitting
shows how this property constrains the rotation frequency.

\begin{lemma}[Splitting]\label{lem:splitting}
{Let} $a\in\nuL$ be \emph{real}, i.e., $a_k\in\R$ for all $k\in\Z$.  Put
\begin{equation}\label{eq:dich_betaRI}
\beta_R\bydef\alpha\cos(\theta-\mu),
\qquad
\beta_I\bydef\alpha\sin(\theta-\mu),
\qquad
d_k\bydef(-1)^k a_k-[g(a)]_k .
\end{equation}
Then $Sa$, $a^\dagger$ and $g(a)$ are real sequences, and for every $k\in\Z$
\begin{equation}\label{eq:dich_splitting}
\operatorname{Re}G_k(a,\mu)=\beta_R\,d_k,
\qquad
\operatorname{Im}G_k(a,\mu)=( {\pi} k+\mu)\,a_k+\beta_I\,d_k .
\end{equation}
\end{lemma}

\begin{proof}
By \eqref{eq:dich_resonance}, $(Sa)_k=e^{-\im {\pi} k}a_k=(-1)^k a_k\in\R$, and
$(a^\dagger)_k= {a_{-k}^\ast}=a_{-k}\in\R$.  A convolution of real sequences is real, so
$g(a)=a\cp a^\dagger\cp Sa$ is real.  Next,
$\lambda_k a_k=\im( {\pi} k+\mu)a_k+\beta(-1)^k a_k$, so that
\[
G_k(a,\mu)=\im( {\pi} k+\mu)a_k+\beta\bigl((-1)^k a_k-[g(a)]_k\bigr)
=\im( {\pi} k+\mu)a_k+(\beta_R+\im\beta_I)\,d_k .
\]
Since $a_k$ and $d_k$ are real, the first summand is purely imaginary.  Taking real and imaginary
parts gives \eqref{eq:dich_splitting}.
\end{proof}

\begin{remark}\label{rem:dich_exact}
The splitting holds in the full sequence space and does not depend on a Fourier
truncation. It uses only the exact resonance \eqref{eq:dich_resonance} and the reality of the
Fourier coefficients.
Here a real sequence means $a_k\in\R$; the profile $w$ itself is generally complex valued.
\end{remark}

The real part of \eqref{eq:dich_splitting} is $\beta_R$ times a quantity independent of
$\mu$. If $\beta_R\ne0$, its vanishing reduces the imaginary equations to a diagonal
system. The following proposition characterizes the resulting exceptional solutions.

\begin{proposition}[Forcing]\label{prop:forcing}
Let $a\in\nuL$ be real and $\mu\in\R$, and suppose $(a,\mu)$ solves the
\emph{full} system $G_k(a,\mu)=0$ for all $k\in\Z$.  Then at least one of the following holds.
\begin{enumerate}[label=\textup{(\roman*)}]
\item $\cos(\theta-\mu)=0$, that is, $\mu\in\theta-\tfrac{\pi}{2}+\pi\Z$.
\item The solution is trivial: either $w\equiv0$, or $w(t)=c\,e^{\im {\pi} k_0 t}$ with
$\abs{c}=1$ and $ {\pi} k_0+\mu=0$ for some $k_0\in\Z$.  In either case
$z(t)=w(t)e^{\im\mu t}$ is an \emph{equilibrium} of \eqref{eq:DDE_complex}, namely $z\equiv0$ or
$z\equiv c$ with $\abs{c}=1$.  Any constant of modulus one is an equilibrium, the factor
$1-\abs{z}^2$ vanishing identically.
\end{enumerate}
\end{proposition}

\begin{proof}
Suppose (i) fails, i.e., $\beta_R=\alpha\cos(\theta-\mu)\ne0$ (recall $\alpha>0$).  By
Lemma~\ref{lem:splitting}, $\operatorname{Re}G_k=\beta_R d_k=0$ forces $d_k=0$ for every
$k\in\Z$, and then $\operatorname{Im}G_k=( {\pi} k+\mu)a_k=0$ for every $k$.  The
linear function $k\mapsto {\pi} k+\mu$ vanishes for at most one integer, and only
if $-\mu/ {\pi}\in\Z$.  We call that integer $k_0$ when it exists.  Hence, $a_k=0$ for all $k\ne k_0$,
so $a=a_{k_0}\delta_{k_0}$ (and $a\equiv0$ if no such integer exists, which is case~(ii) with
$w\equiv0$).

For a one-mode real sequence $a=c\,\delta_{k_0}$, $c\in\R$, we have
$a^\dagger=c\,\delta_{-k_0}$ and $Sa=(-1)^{k_0}c\,\delta_{k_0}$, so
$a\cp a^\dagger=c^2\delta_0$ and
\[
g(a)=c^2\delta_0\cp\bigl((-1)^{k_0}c\,\delta_{k_0}\bigr)=(-1)^{k_0}c^3\,\delta_{k_0},
\qquad\text{i.e.}\qquad
[g(a)]_{k_0}=(-1)^{k_0}c^3 .
\]
Therefore $d_{k_0}=(-1)^{k_0}c-(-1)^{k_0}c^3=(-1)^{k_0}c\,(1-c^2)$, and $d_{k_0}=0$ gives
$c\in\{0,\pm1\}$.  If $c=0$, then $w\equiv0$.  If $\abs{c}=1$, then
$w(t)=c\,e^{\im {\pi} k_0t}$ and $z(t)=w(t)e^{\im\mu t}=c\,e^{\im( {\pi} k_0+\mu)t}\equiv c$.
This is case~(ii).
\end{proof}

\begin{corollary}[Allowed rotation frequencies]\label{cor:dich_forced_mu}
{Let} $(a,\mu)$ solve the full system with real Fourier coefficients and
suppose that $z(t)=w(t)e^{\im\mu t}$ is not an equilibrium. Then
$\mu=\theta-\pi/2+j\pi$ for some $j\in\Z$.
\end{corollary}

Let $\mu=\theta-\pi/2+j\pi$ be a frequency from
Corollary~\ref{cor:dich_forced_mu}. With $\mu_0=\theta-\pi/2$, set
\[
\widehat w(t)=e^{\im j\pi t}w(t).
\]
Then $z(t)=e^{\im\mu_0t}\widehat w(t)$ represents the same solution of the original
DDE. The new profile is still $2$-periodic and has real Fourier coefficients
${\widehat a}_k=a_{k-j}$, with
$\norm{{\widehat a}}_\nu\le\nu^{|j|}\norm a_\nu$.
Thus every nonequilibrium solution in this symmetric class admits the representative
\begin{equation}\label{eq:dich_forced_mu}
\mu=\theta-\frac\pi2,
\qquad \beta=\im\alpha,
\qquad \beta_R=0,\quad\beta_I=\alpha.
\end{equation}
We use this representative throughout the validation. The nonvanishing coefficient
conditions below are checked in this frame. The value of $\mu$ equals the angular
frequency of the circular orbit in Theorem~\ref{thm:P2}, but the profile $w$ need not
be stationary.

The conjugation-reversal underlying the real coefficient restriction also preserves
the periodic equation at this frequency.

\begin{lemma}[Reversibility]\label{lem:reversibility}
On $2$-periodic functions, define the conjugation-reversal
\begin{equation}\label{eq:dich_rho}
({\sigma} w)(t)\bydef {w(-t)^\ast},
\qquad\text{acting on Fourier coefficients by}\qquad
({\sigma} a)_k= {a_k^\ast}.
\end{equation}
If $w$ is a $2$-periodic solution of the $w$-DDE \eqref{eq:w_dde} with coupling coefficient $\beta$,
then $\sigma w$ is a $2$-periodic solution of the same equation with $\beta$ replaced by
$- {\beta^\ast}$.  Consequently, $\sigma$ is a symmetry of the $2$-periodic problem if and only if
$\beta=- {\beta^\ast}$, i.e., $\beta\in\im\R$, i.e., $\cos(\theta-\mu)=0$, the frequency condition in Proposition~\ref{prop:forcing}.  The fixed-point set of $\sigma$ is the real
subspace $\{a:\ a_k\in\R\ \forall k\}$, i.e., the profiles with $w(-t)= {w(t)^\ast}$.
\end{lemma}

\begin{proof}
Write \eqref{eq:w_dde} as $\dot w(t)=-\beta\bigl(1-\abs{w(t)}^2\bigr)w(t-1)-\im\mu\,w(t)$ and set
$v\bydef\sigma w$, so that $v(t)= {w(-t)^\ast}$ and $\dot v(t)=- {\dot w(-t)^\ast}$.  Then
\[
\dot v(t)
= {\Bigl(\beta\bigl(1-\abs{w(-t)}^2\bigr)w(-t-1)+\im\mu\,w(-t)\Bigr)^\ast}
= {\beta^\ast}\,\bigl(1-\abs{v(t)}^2\bigr)\, {w(-t-1)^\ast}-\im\mu\,v(t),
\]
using $\abs{w(-t)}=\abs{v(t)}$.  Now $ {w(-t-1)^\ast}= {w(-(t+1))^\ast}=v(t+1)$, and
$2$-periodicity turns the \emph{advance} into the \emph{delay},
\[
v(t+1)=v(t+1-2)=v(t-1).
\]
Hence, $\dot v(t)=-\bigl(- {\beta^\ast}\bigr)\bigl(1-\abs{v(t)}^2\bigr)v(t-1)-\im\mu\,v(t)$,
which is \eqref{eq:w_dde} with $- {\beta^\ast}$ in place of $\beta$.  The remaining assertions are
immediate: $\beta=- {\beta^\ast}$ says $\operatorname{Re}\beta=\alpha\cos(\theta-\mu)=0$, and
$(\sigma a)_k= {a_k^\ast}$ follows by conjugating
$w(-t)=\sum_k a_k e^{-\im {\pi} kt}$ term by term.
\end{proof}

\begin{remark}\label{rem:dich_reversibility}
A delay normally destroys time-reversal symmetry, since reversing time turns a delay into an
advance.  Here the symmetry is restored, but only on the class of $2$-periodic profiles (where the
advance folds back onto the delay) and only for purely imaginary $\beta$.  We restrict the existence proof to $\sigma$-symmetric profiles.  In the rotating frame, $\sigma$-symmetry says exactly
that $\operatorname{Re}w$ is even and $\operatorname{Im}w$ is odd in $t$.  Section~\ref{sec:cap} validates solutions in this symmetric class.
\end{remark}

\subsection{The reduced real system}\label{subsec:dich_real}

With the exact choice \eqref{eq:dich_forced_mu}, Lemma~\ref{lem:splitting} gives
$G_k(a,\mu)=\im F^R_k(a)$ on $X$, where
\begin{equation}\label{eq:dich_FR}
F^R_k(a)\bydef\lambda_k^R a_k-\alpha[g(a)]_k,
\qquad \lambda_k^R\bydef\pi k+\mu+\alpha(-1)^k\in\R,
\qquad k\in\Z.
\end{equation}
Consequently, solving $F^R(a)=0$ on the real space $X$ satisfies every real and
imaginary equation of the full system. The regularity argument in
Section~\ref{subsec:dde_zero} then gives a classical DDE solution.

Truncating to $|k|\le N$ gives $2N+1$ real equations for $2N+1$ real unknowns.
The continuous rotations and time translations of the complex problem reduce to
discrete symmetries on this subspace: a rotation preserves a nonzero real sequence
only for angles in $\pi\Z$, and, when $a_1\ne0$, a time translation preserves reality
only for shifts in $\Z$. The real system therefore needs no additional phase conditions.
Its local nonsingularity is established by the interval estimates in
Section~\ref{sec:cap}.

\subsection{Star-shaped periodic solutions and orbit closure}\label{subsec:dich_geometry}

We fix the branch \eqref{eq:dich_forced_mu} and determine its rotation number.  Let
$\theta=\pi/n$ with $n>2$.  When $n\in\mathbb{Q}$, define coprime positive integers $p,m$ by
\begin{equation}\label{eq:dich_star_numbers}
\frac{m}{p}=\frac{n-2}{2n}\quad\text{in lowest terms};
\qquad\text{for integer }n,\quad
p=\frac{2n}{\gcd(2n,\,n-2)},
\quad
m=\frac{n-2}{\gcd(2n,\,n-2)}.
\end{equation}
Then $\gcd(p,m)=1$, and $p>2m$ because $p/m=2n/(n-2)>2$ for $n>2$.  On this branch
$\mu=\theta-\pi/2$ gives, exactly,
\begin{equation}\label{eq:dich_rotation_number}
\frac{\mu}{\pi}=\frac{1}{n}-\frac12=\frac{2-n}{2n},
\qquad\text{which equals }-\frac{m}{p}\in\mathbb{Q}\ \text{ when }n\in\mathbb{Q},
\end{equation}
and therefore, for $n\in\mathbb{Q}$, by the relative periodicity \eqref{eq:relative_periodicity},
\begin{equation}\label{eq:dich_multiplier}
z(t+2)=e^{-2\pi\im m/p}\,z(t):
\end{equation}
translation by one profile period rotates the solution by $-2\pi m/p$, which is the symmetry
\eqref{eq:star_shaped_symmetry} in Definition~\ref{def:star_shaped}.

\begin{proposition}[Star-shaped periodic solutions and orbit closure]\label{prop:dich_star}
Let $\theta=\pi/n$ with $n>2$ and let $\tw$ be a $2$-periodic solution of \eqref{eq:w_dde} on the
branch $\mu=\theta-\pi/2$ whose Fourier coefficients $\ta$ satisfy $\ta_0\ne0$ and $\ta_1\ne0$.
Put $z(t)=\tw(t)e^{\im\mu t}$. Then $|z(t)|$ is nonconstant, and the following hold.
\begin{enumerate}[label=\textup{(\alph*)}]
\item \textup{(rational case)} If $n\in\mathbb{Q}$, with $p,m$ as in
\eqref{eq:dich_star_numbers}, then $z$ is a star-shaped periodic solution of type $(p,m)$
in the sense of Definition~\ref{def:star_shaped}. Its fundamental period is $L\bydef 2p$.
\item \textup{(Fourier sparsity)} In the same situation, the Fourier support of $z$ is contained in the frequencies
$\Omega(pk-m)$, $k\in\Z$, where $\Omega\bydef2\pi/L$. Thus its signed indices lie in
$-m\pmod p$. The real coordinates $x$ and $y$ have support contained in the two
residue classes $\pm m\pmod p$. The possible positive index magnitudes are
$\{m,p-m,p+m,2p-m,\dots\}$.
\item \textup{(irrational case)} If $n\notin\mathbb{Q}$, then $\mu/\pi=(2-n)/(2n)$ is
irrational, $z$ never closes, and the closure of the single orbit $\{z(t):t\in\R\}$ is the annulus
\[
\Bigl\{\varrho\,e^{\im\varphi}:\ \varrho\in\bigl[\min_t\abs{\tw(t)},\ \max_t\abs{\tw(t)}\bigr],\
\varphi\in[0,2\pi)\Bigr\}.
\]
\end{enumerate}
\end{proposition}

\begin{proof}
We first show that the modulus is nonconstant. If $|\tw(t)|\equiv\varrho$, then
$g(\ta)=\varrho^2S\ta$. Since $\beta=\im\alpha$ and $\ta_0, \ta_1\ne0$,
the $k=0,1$ equations give
\[
\mu+\alpha(1-\varrho^2)=0,
\qquad \pi+\mu-\alpha(1-\varrho^2)=0.
\]
Adding them yields $\pi+2\mu=2\theta=0$, contradicting $\theta=\pi/n>0$.
Thus $|z|=|\tw|$ is nonconstant. For rational $n$, this fact and
\eqref{eq:dich_multiplier} give Definition~\ref{def:star_shaped}.

Assume now $n\in\mathbb{Q}$. Expanding $\tw(t)=\sum_k \ta_k e^{\im {\pi} kt}$ gives
$z(t)=\sum_k \ta_k e^{\im( {\pi} k+\mu)t}$, and by \eqref{eq:dich_rotation_number}
and $\Omega=2\pi/L=\pi/p$,
\[
 {\pi} k+\mu=\pi k-\frac{m\pi}{p}=\Omega\,(pk-m),
\qquad k\in\Z .
\]
Thus the $z$-frequencies have indices in $-m$ modulo $p$. Since
$x=(z+ {z^\ast})/2$ and $y=(z- {z^\ast})/(2\im)$, their supports are contained in the
union of that class and its negative. This proves (b).  In particular, $z(t+L)=z(t)$, so $L$ is a
period.

For (a), let $P>0$ be any period of $z$. Uniqueness of Fourier coefficients and
$\ta_0, \ta_1\ne0$ imply
\[
e^{\im\mu P}=e^{\im(\pi+\mu)P}=1.
\]
Dividing these identities gives $e^{\im\pi P}=1$, so $P=2q$ for a positive integer $q$.
Since $\mu=-m\pi/p$, the first identity implies $mq/p\in\Z$. Coprimality gives
$p\mid q$, proving that $L=2p$ is the fundamental period. Finally,
\eqref{eq:dich_multiplier} makes the orbit invariant under rotation by $-2\pi m/p$.
This rotation generates the same cyclic group as rotation by $2\pi/p$.

For (c), $(2-n)/(2n)=1/n-1/2$ is irrational whenever $n$ is, so $\mu/\pi\notin\mathbb{Q}$ and
$\{e^{2\im\mu\ell}\}_{\ell\in\Z}$ is dense in the unit circle.  By
\eqref{eq:relative_periodicity}, $z(t+2\ell)=e^{2\im\mu\ell}z(t)$, so for each fixed $t$ the
closure of $\{z(t+2\ell):\ell\in\Z\}$ is the whole circle of radius $\abs{\tw(t)}$.  As $t$ ranges
over $[0,2]$, the continuous $2$-periodic function $\abs{\tw}$ attains exactly the values in
$[\min_t\abs{\tw(t)},\max_t\abs{\tw}(t)]$.  Taking the union of the corresponding circles gives the
stated annulus, and $\abs{z(t)}=\abs{\tw(t)}$ shows the orbit is contained in it.  Finally $z$ is
not periodic: since $\ta_0\ne0$ and $\ta_1\ne0$, the frequencies $\mu$ and $ {\pi+\mu}$
both occur in the spectrum of $z$, whereas a $P$-periodic function has spectrum contained in
$(2\pi/P)\Z$.  That would force
$(\pi+\mu)/\mu=1+2n/(2-n)=(2+n)/(2-n)\in\mathbb{Q}$, hence $n\in\mathbb{Q}$, a contradiction.
\end{proof}

\begin{remark}\label{rem:dich_closure}
The rotation identity \eqref{eq:dich_rotation_number} gives
\[
p(\pi-2\pi/n)=p(-2\mu)=2\pi m.
\]  The hypotheses $\ta_0\ne0$, $\ta_1\ne0$ of
Proposition~\ref{prop:dich_star} are \emph{open} conditions, hence verifiable a posteriori from a
validated enclosure: if \akt{an approximate solution $\ba$ and the radii-polynomial radius $r$
satisfy $\norm{\ta-\ba}_\nu\le r$ and}
$\abs{\ba_k}>r/\nu^{\abs{k}}$ for $k=0,1$, then $\ta_k\ne0$ for $k=0,1$.
\end{remark}

\begin{table}[htbp]\centering
\caption{Rotation data on the branch $\mu=\theta-\pi/2$, $\theta=\pi/n$: $m/p=(n-2)/(2n)$ in lowest
terms, fundamental period $L=2p$ under the nonvanishing hypotheses of
Proposition~\ref{prop:dich_star}, and allowed Fourier indices $pl\pm m$, $l\in\Z$.
The entries follow from the exact rotation number.}
\label{tab:dich_starnumbers}
\begin{tabular}{cccccc}
\hline
$n$ & $\theta=\pi/n$ & $p$ & $m$ & $L=2p$ & allowed indices \\
\hline
$3$  & $\pi/3$  & $6$  & $1$  & $12$  & $6l\pm1$ \\
$4$  & $\pi/4$  & $4$  & $1$  & $8$   & $4l\pm1$ \\
$5$  & $\pi/5$  & $10$ & $3$  & $20$  & $10l\pm3$ \\
$6$  & $\pi/6$  & $3$  & $1$  & $6$   & $3l\pm1$ \\
$7$  & $\pi/7$  & $14$ & $5$  & $28$  & $14l\pm5$ \\
$10$ & $\pi/10$ & $5$  & $2$  & $10$  & $5l\pm2$ \\
$35$ & $\pi/35$ & $70$ & $33$ & $140$ & $70l\pm33$ \\
\hline
\end{tabular}
\end{table}

\begin{remark}[Winding number]\label{rem:dich_table_correction}
If $w$ avoids the origin, then over $L=2p$ the signed winding number satisfies
\[
\operatorname{wind}(z;0)=p\,\operatorname{wind}(w;0)-m,
\]
where the winding of $w$ is taken over one profile period.
In particular, if $\operatorname{Re}w(t)>0$ for every $t$, then $w$ is homotopic to a
positive constant within the right half-plane and $\operatorname{wind}(z;0)=-m$.
We verify this inequality for the seven validated profiles below. Thus their winding
number is determined in addition to the symmetry in Definition~\ref{def:star_shaped}.
\end{remark}

\section{Computer-assisted proof}\label{sec:cap}

We now validate zeros of the real system \eqref{eq:dich_FR}. For each parameter pair,
we compute an approximate profile, bound the residual and derivatives in $X$ defined in \eqref{eq:real_space}, and
apply the Radii Polynomial Theorem. The estimates include all Fourier modes. A single
validation then provides the profile needed for the periodicity and winding arguments
of Section~\ref{subsec:dich_geometry}.

Throughout this section $\alpha$, $\theta$ and $\mu=\theta-\pi/2$ in \eqref{eq:dich_forced_mu} are fixed data.
For the seven examples, $\alpha$ is rational and $\theta$ and $\mu$ are rational
multiples of $\pi$. The computation retains these rational inputs and uses a
rigorous enclosure of $\pi$. The identities $S=\operatorname{diag}((-1)^k)$ in \eqref{eq:shift_dagger}
and $\beta=\im\alpha$ are imposed explicitly.

\subsection{The validation map}\label{subsec:dde_rpa}\label{subsec:dde_precond}
Let $N\ge1$, and let $\Pi_N$ and $\Pi_{>N}=I-\Pi_N$ project onto the modes
$|k|\le N$ and $|k|>N$, respectively. We choose $N$ so that
\begin{equation}\label{eq:N_choice}
\pi(N+1)>|\mu|+\alpha.
\end{equation}
For all $k\in\Z$, the diagonal $\lambda_k^R$ in \eqref{eq:dich_FR} satisfies
$|\lambda_k^R|\ge\pi|k|-|\mu|-\alpha$. Hence \eqref{eq:N_choice} gives $\lambda_k^R\ne0$ for
all $|k|>N$, and since $|\lambda_k^R|\to\infty$ as $|k|\to\infty$, we also have
$\inf_{|k|>N}|\lambda_k^R|>0$. For the seven cases treated in
Section~\ref{subsec:cap_results}, $N=30$ satisfies \eqref{eq:N_choice}.

\paragraph{The map $F$.}
Since $F^R$ is unbounded on $X$, we precondition its tail to define a
smooth map $F:X\to X$, with $g$ in \eqref{eq:g_cubic}:
\begin{equation}\label{eq:F_precond}
F_k(a)\bydef
\begin{cases}
F^R_k(a),& |k|\le N,\\[2pt]
a_k-\dfrac{\alpha}{\lambda_k^R}[g(a)]_k,& |k|>N.
\end{cases}
\end{equation}
The finite part is polynomial. In the tail, the diagonal operator
$\operatorname{diag}(\alpha/\lambda_k^R)_{|k|>N}$ is bounded on $X$ because
$\inf_{|k|>N}|\lambda_k^R|>0$. Hence $F$ is $C^\infty$ on $X$.
The tail diagonal is nonzero, so $F(a)=0$ if and only if all the equations
$F^R_k(a)=0$ hold. Every such zero defines a classical solution by
Section~\ref{subsec:dde_zero}.

\paragraph{The approximate zero and the operator $A$.}
Numerical integration of \eqref{eq:w_dde}
and a discrete Fourier transform give an initial guess. We take its real coefficient
part and apply Newton's method to the truncation of \eqref{eq:dich_FR} to $|k|\le N$, obtaining
$\ba\in\Pi_NX$. This finite sequence is embedded in $X$ by padding zeros to all higher modes. The subsequent validation bounds its error, including the infinite tail.

Write
\[
J_N\bydef\Pi_NDF^R(\ba)\Pi_N
\]
for the finite Jacobian of $F^R$ defined in \eqref{eq:dich_FR}, regarded as a $(2N+1)\times(2N+1)$ matrix acting on $\Pi_NX$, and let $B$ be a floating-point approximate inverse of the
midpoint of an interval enclosure of $J_N$. Its entries are treated as fixed, exactly
representable numbers. We verify its accuracy below, rather than assuming that the
numerical inversion is exact.

We use the bounded approximate inverse
\begin{equation}\label{eq:A_split}
A\bydef B\Pi_N\oplus\Pi_{>N}:X\to X.
\end{equation}
The finite estimate $\norm{I-BJ_N}<1$ below proves that $B$ is nonsingular,
and hence that $A$ is injective. We apply the following theorem to the Newton-like map
$a\mapsto a-AF(a)$.

\subsection{The Radii Polynomial Theorem}\label{subsec:dde_thm}
\begin{theorem}[Radii Polynomial Theorem]\label{thm:rpa}
Let $X$ and $\mathcal Y$ be real Banach spaces, let $F:X\to\mathcal Y$ be $C^2$,
let $\bx\in X$ and let $A:\mathcal Y\to X$ be an injective bounded linear operator.
Write $B(\bx,r)$ for the closed ball in $X$.  Fix $R>0$ and suppose that nonnegative numbers $Y$, $Z_0$,
$Z_1$ satisfy
\begin{equation}\label{eq:YZ}
\norm{AF(\bx)}_X\le Y,
\qquad
\norm{I-A\,DF(\bx)}_{\mathscr B(X)}\le Z_0,
\qquad
\sup_{x\in B(\bx,R)}\norm{A\,D^2F(x)}_{\mathscr B^2(X)}\le Z_1 .
\end{equation}
Define the radii polynomial
\begin{equation}\label{eq:radiipoly}
p(r)\bydef Y+(Z_0-1)\,r+\frac{Z_1}{2}\,r^2 .
\end{equation}
If $0\le r\le R$, $p(r)\le0$ and $Z_0+Z_1r<1$, then $F$ has a unique zero $\tx$ in
$B(\bx,r)$.
\end{theorem}

\begin{proof}
Set $T(x)\bydef x-AF(x)$. Taylor's formula and \eqref{eq:YZ} give
$\norm{T(x)-\bx}_X\le Y+Z_0r+\tfrac12Z_1r^2\le r$ on $B(\bx,r)$.
Moreover, $\norm{DT(x)}\le Z_0+Z_1r<1$ there. The contraction mapping theorem gives a
unique fixed point, which is a zero of $F$ because $A$ is injective.
\end{proof}

\begin{remark}[Admissible radius]\label{rem:admissible_radius}
We now describe how an admissible radius in Theorem~\ref{thm:rpa} is found. Let $Z_1>0$, $Z_0<1$ and
$\Delta\bydef(1-Z_0)^2-2Z_1Y>0$. Then $p(r)$ in \eqref{eq:radiipoly} has two real roots
\[
r_\pm=\frac{(1-Z_0)\pm\sqrt{\Delta}}{Z_1},
\]
and the set of $r\ge0$ with $p(r)\le0$ is the interval $[r_-,r_+]$. The smaller root gives the
tightest enclosure, and the contraction condition holds there automatically, since
$Z_1r_-=(1-Z_0)-\sqrt{\Delta}$ implies $Z_0+Z_1r_-=1-\sqrt{\Delta}<1$. Hence, for $Z_1>0$,
an admissible radius exists if and only if $Z_0<1$, $\Delta>0$ and $r_-\le R$, and we take
$r=r_-$. For the computation we use the equivalent form
\[
r_-=\frac{2Y}{1-Z_0+\sqrt{\Delta}},
\]
which is obtained by rationalizing the numerator and avoids cancellation in the smaller quadratic root. If $\Delta=0$, the
only radius with $p(r)\le0$ has contraction bound equal to one and is inadmissible.
For $Z_1=0$, the conditions are $Z_0<1$ and $Y/(1-Z_0)\le R$.
\end{remark}

For the application here, both spaces in Theorem~\ref{thm:rpa} are the real sequence
space \eqref{eq:real_space}, and the center is $\ba$. Local uniqueness is asserted
within this real space and the validated ball.

\subsection{The bounds \texorpdfstring{$Y$, $Z_0$ and $Z_1$}{Y, Z0 and Z1}}\label{subsec:dde_bounds}

We give the estimates for the map \eqref{eq:F_precond} and the operator
\eqref{eq:A_split}. All finite vectors and matrices are enclosed with outward-rounded
interval arithmetic. The sequence norm is the weighted $\ell^1$ norm in \eqref{eq:ell1nu_complex}; for a real
matrix $M$ its induced norm is the maximum weighted column sum
\[
\norm M=\max_j\frac{1}{\nu^{|j|}}\sum_k\nu^{|k|}|M_{kj}|.
\]

\paragraph{A lower bound on the tail diagonal.}
For all $k\in\Z$, we have $|\lambda_k^R|\ge\pi|k|-|\mu|-\alpha$. By \eqref{eq:N_choice}, the
right-hand side is positive for $|k|>N$, and hence the tail diagonal does not vanish. To obtain a
sharper lower bound, we fix an integer $K_\ast\ge N+1$ and use
$|\lambda_k^R|\ge\pi(K_\ast+1)-|\mu|-\alpha>0$ for $|k|>K_\ast$. Hence a number
$\lambda_{\min}$ computed by a finite interval scan of $N<|k|\le K_\ast$ satisfies
\begin{equation}\label{eq:lambda_min}
\lambda_{\min}\le
\min\left\{\min_{N<|k|\le K_\ast}|\lambda_k^R|,
\ \pi(K_\ast+1)-|\mu|-\alpha\right\}
\le\inf_{|k|>N}|\lambda_k^R|.
\end{equation}
Under \eqref{eq:N_choice}, the scan gives $\lambda_{\min}>0$. For the seven cases treated
in Section~\ref{subsec:cap_results}, we obtain $\lambda_{\min}>90$. This check covers the whole tail. It does not require division by
any finite-mode $\lambda_k^R$.
In the tail, the preconditioner
$\Lambda_{>N}^{-1}$, where $\Lambda_{>N}\bydef\operatorname{diag}(\lambda_k^R)_{|k|>N}$, has norm at most $1/\lambda_{\min}$, and hence the
diagonal operator $\operatorname{diag}(\alpha/\lambda_k^R)_{|k|>N}$ multiplying $g(a)$ has norm
at most $\alpha/\lambda_{\min}$. The latter bound enters $C_{TT}$ in \eqref{eq:Z0bound}, and the
former enters $Z_1$ in \eqref{eq:Z1bound}.

\paragraph{The residual bound $Y$.}
The norm on $X$ splits as $\norm u_\nu=\norm{\Pi_Nu}_\nu+\norm{\Pi_{>N}u}_\nu$. By
\eqref{eq:A_split}, $A$ acts by $B$ on the finite modes and by the identity on the tail, and
by \eqref{eq:F_precond}, $\Pi_NF(\ba)=\Pi_NF^R(\ba)$. Hence
\begin{align*}
\norm{AF(\ba)}_\nu
&=\norm{B\Pi_NF(\ba)}_\nu+\norm{\Pi_{>N}F(\ba)}_\nu\\
&=\norm{B\Pi_NF^R(\ba)}_\nu
+\sum_{|k|>N}\nu^{|k|}\Bigl|\ba_k-\frac{\alpha}{\lambda_k^R}[g(\ba)]_k\Bigr|\\
&=\norm{B\Pi_NF^R(\ba)}_\nu
+\sum_{N<|k|\le3N}\nu^{|k|}\frac{\alpha|[g(\ba)]_k|}{|\lambda_k^R|}.
\end{align*}
In the last step, we used that $\ba_k=0$ for $|k|>N$ and that the cubic sequence $g(\ba)$ is supported on $|k|\le3N$. The right-hand side is a finite sum, and we take $Y$ as an upper
bound of it computed with interval arithmetic:
\begin{equation}\label{eq:Ybound}
\norm{AF(\ba)}_\nu=\norm{B\Pi_NF^R(\ba)}_\nu+
\sum_{N<|k|\le3N}\nu^{|k|}\frac{\alpha|[g(\ba)]_k|}{|\lambda_k^R|}
\le Y.
\end{equation}
Thus the entire residual, including its tail, is obtained from a finite computation.

\paragraph{The derivative bound $Z_0$.}
To control $M\bydef I-ADF(\ba)$, we separate finite and tail input columns. With $S$ and $\dagger$ in \eqref{eq:shift_dagger}, let
\[
Q_1\bydef\ba^\dagger\cp S\ba,\qquad
Q_2\bydef\ba\cp S\ba,\qquad
Q_3\bydef\ba\cp\ba^\dagger.
\]
Each kernel is supported on $|k|\le2N$. The derivative $K\bydef Dg(\ba)$ of \eqref{eq:g_cubic} acts by
\begin{equation}\label{eq:dich_Dg}
Kh=h\cp Q_1+h^\dagger\cp Q_2+(Sh)\cp Q_3,
\qquad
K_{kj}=(Q_1)_{k-j}+(Q_2)_{k+j}+(-1)^j(Q_3)_{k-j}.
\end{equation}
In particular,
$\norm K\le\norm{Q_1}_\nu+\norm{Q_2}_\nu+\norm{Q_3}_\nu$.
For $|k|,|j|\le N$, \eqref{eq:dich_FR} gives $(J_N)_{kj}=\lambda_k^R\delta_{kj}-\alpha K_{kj}$.
More generally, \eqref{eq:F_precond} gives the rows of $DF(\ba)$ for all columns $j\in\Z$:
\[
(DF(\ba))_{kj}=
\begin{cases}
\lambda_k^R\delta_{kj}-\alpha K_{kj}, & |k|\le N,\\[2pt]
\delta_{kj}-\dfrac{\alpha}{\lambda_k^R}K_{kj}, & |k|>N.
\end{cases}
\]
Since $A$ acts by $B$ on the finite rows and by the identity on the tail rows, the identity
cancels in the tail rows of $M=I-ADF(\ba)$, and only the nonlinear part remains there. In the
finite rows, the columns $|j|\le N$ give $I-BJ_N$, while a tail column $|j|>N$ has
$\delta_{kj}=0$ for $|k|\le N$ and gives $B\Pi_N\alpha Ke_j$.
With $\Lambda_{>N}$ as above, $M$ has the block form
\begin{equation}\label{eq:M_block}
M=
\begin{pmatrix}
M_{FF} & M_{FT}\\
M_{TF} & M_{TT}
\end{pmatrix}
=
\begin{pmatrix}
I-BJ_N & \alpha B\Pi_NK\Pi_{>N}\\[2pt]
\alpha\Lambda_{>N}^{-1}\Pi_{>N}K\Pi_N & \alpha\Lambda_{>N}^{-1}\Pi_{>N}K\Pi_{>N}
\end{pmatrix},
\end{equation}
where the first block row and column correspond to the modes $|k|\le N$ and the second to
the modes $|k|>N$.
For finite input $|j|\le N$, the tail output vanishes when $|k|>3N$.
For tail input $|j|>3N$, the finite output vanishes. Both statements follow from
\eqref{eq:dich_Dg}, since $K_{kj}\ne0$ requires $|k-j|\le2N$ or $|k+j|\le2N$. We can therefore compute the
finite-column bound and the finite output from tail columns by finite sums. With
$e_j$ denoting the unit coefficient at mode $j$, set
\begin{equation}\label{eq:Z0bound}
\begin{aligned}
C_{\mathrm{fin}}&\bydef\max_{|j|\le N}\frac{1}{\nu^{|j|}}
\left(\sum_{|k|\le N}\nu^{|k|}|(M_{FF})_{kj}|
+\sum_{N<|k|\le3N}\nu^{|k|}\frac{\alpha|K_{kj}|}{|\lambda_k^R|}\right),\\
C_{FT}&\bydef\max_{N<|j|\le3N}
\frac{\norm{B\Pi_N\alpha K e_j}_\nu}{\nu^{|j|}},\\
C_{TT}&\bydef\frac{\alpha}{\lambda_{\min}}
\left(\norm{Q_1}_\nu+\norm{Q_2}_\nu+\norm{Q_3}_\nu\right).
\end{aligned}
\end{equation}
We now derive the bound on $M=I-ADF(\ba)$. The operator norm on $X$ is the maximum
weighted column sum, and we split each column sum into its finite rows $|k|\le N$ and its tail
rows $|k|>N$. Writing
\[
c_j\bydef\frac{1}{\nu^{|j|}}\sum_{k\in\Z}\nu^{|k|}|M_{kj}|,
\qquad
\norm{I-ADF(\ba)}=\sup_{j\in\Z}c_j,
\]
we have, for a finite column $|j|\le N$, by \eqref{eq:M_block} and the support of $K$,
\[
c_j=\frac{1}{\nu^{|j|}}\left(\sum_{|k|\le N}\nu^{|k|}|(M_{FF})_{kj}|
+\sum_{N<|k|\le3N}\nu^{|k|}\frac{\alpha|K_{kj}|}{|\lambda_k^R|}\right)
\le C_{\mathrm{fin}},
\]
where the two sums are the column sums of $M_{FF}$ and $M_{TF}$.
For a tail column $|j|>N$, the two terms below are the column sums of $M_{FT}$ and $M_{TT}$.
The first vanishes unless $|j|\le3N$, and $\nu\ge1$ gives
$\nu^{|k|}\le\nu^{|k-j|}\nu^{|j|}$ and $\nu^{|k|}\le\nu^{|k+j|}\nu^{|j|}$. Hence
\eqref{eq:lambda_min} and \eqref{eq:dich_Dg} yield
\begin{align*}
c_j
&=\frac{\norm{B\Pi_N\alpha Ke_j}_\nu}{\nu^{|j|}}
+\frac{1}{\nu^{|j|}}\sum_{|k|>N}\nu^{|k|}\frac{\alpha|K_{kj}|}{|\lambda_k^R|}\\
&\le C_{FT}+\frac{\alpha}{\lambda_{\min}}\sum_{k\in\Z}\Bigl(\nu^{|k-j|}|(Q_1)_{k-j}|
+\nu^{|k+j|}|(Q_2)_{k+j}|+\nu^{|k-j|}|(Q_3)_{k-j}|\Bigr)\\
&=C_{FT}+C_{TT}.
\end{align*}
Taking the supremum over $j$ and enclosing $C_{\mathrm{fin}}$, $C_{FT}$ and $C_{TT}$ with interval
arithmetic, we obtain
\begin{equation}\label{eq:Z0combine}
Z_0\bydef\max\{C_{\mathrm{fin}},C_{FT}+C_{TT}\}
\ \ge\ \norm{I-ADF(\ba)}.
\end{equation}
Since
$\norm{M_{FF}}\le C_{\mathrm{fin}}\le Z_0$, the check $Z_0<1$ also establishes the
injectivity needed in Theorem~\ref{thm:rpa}.

\paragraph{The second derivative bound $Z_1$.}
Only the cubic nonlinearity contributes to the second derivative, since $\mu$ is fixed.
The map $g(a)=a\cp a^\dagger\cp Sa$ in \eqref{eq:g_cubic} is the restriction to the diagonal of the real
trilinear map $(x,y,z)\mapsto x\cp y^\dagger\cp Sz$. Hence $D^2g(a)[h,v]$ is the sum of the
six terms obtained by placing $h$ and $v$ in two of the three slots and $a$ in the remaining
one, for example $h\cp v^\dagger\cp Sa$ and $a\cp h^\dagger\cp Sv$. Since $\dagger$ and $S$
are isometries of $\nuL$, the Banach-algebra property gives, for $a\in B(\ba,R)$,
\[
\norm{D^2g(a)[h,v]}_\nu\le6\norm a_\nu\norm h_\nu\norm v_\nu
\le6(\norm\ba_\nu+R)\norm h_\nu\norm v_\nu.
\]
For $a\in B(\ba,R)$ and $h,v\in X$, write $z\bydef D^2g(a)[h,v]$. By \eqref{eq:F_precond},
the finite part of $D^2F(a)[h,v]$ is $-\alpha\Pi_Nz$ and its tail part is
$-\alpha\Lambda_{>N}^{-1}\Pi_{>N}z$, and $A$ in \eqref{eq:A_split} applies $B$ to the finite
part. Splitting the norm into the finite and tail parts, we obtain
\begin{align*}
\norm{AD^2F(a)[h,v]}_\nu
&=\alpha\norm{B\Pi_Nz}_\nu+\alpha\norm{\Lambda_{>N}^{-1}\Pi_{>N}z}_\nu\\
&\le\alpha\norm{B}\norm{\Pi_Nz}_\nu+\alpha\lambda_{\min}^{-1}\norm{\Pi_{>N}z}_\nu\\
&\le\alpha\max\{\norm{B},\lambda_{\min}^{-1}\}\norm z_\nu\\
&\le6\alpha(\norm\ba_\nu+R)\max\{\norm{B},\lambda_{\min}^{-1}\}\norm h_\nu\norm v_\nu.
\end{align*}
The second line uses \eqref{eq:lambda_min}, the third uses
$\norm z_\nu=\norm{\Pi_Nz}_\nu+\norm{\Pi_{>N}z}_\nu$, and the last uses the bound on $D^2g(a)$
above. Taking the supremum over $\norm h_\nu,\norm v_\nu\le1$ and $a\in B(\ba,R)$, and enclosing
the right-hand side with interval arithmetic, we obtain
\begin{equation}\label{eq:Z1bound}
Z_1\bydef6\alpha(\norm\ba_\nu+R)
\max\{\norm{B},\lambda_{\min}^{-1}\}
\ \ge\ \sup_{a\in B(\ba,R)}\norm{AD^2F(a)}.
\end{equation}

\subsection{Validated profiles and proof of Theorem~\ref{thm:main}}\label{subsec:cap_results}

For each of the seven parameter pairs in Theorem~\ref{thm:main}, we use $N=30$,
$\nu=1$ and $R=10^{-2}$ in Theorem~\ref{thm:rpa}. The bounds below satisfy the hypotheses of
Theorem~\ref{thm:rpa} for an admissible radius $r$. Hence each ball $B(\ba,r)\subset X$
contains a unique zero $\ta$ of $F$. By construction and regularity (Section~\ref{subsec:dde_zero}), its Fourier
series $\tw$ solves the full rotating-frame equation \eqref{eq:w_dde} at the exact nominal parameters.

\begin{table}[htbp]\centering
\caption{Validated bounds for the real system \eqref{eq:dich_FR} with $\mu=\theta-\pi/2$, $N=30$, $\nu=1$ and
$R=10^{-2}$ in Theorem~\ref{thm:rpa}. The values of $Y$, $Z_0$ and $Z_1$ are upper bounds, and $r$ is an admissible radius:
both $p(r)<0$ (see \eqref{eq:radiipoly}) and $Z_0+Z_1r<1$ hold with the tabulated values.}\vspace*{1mm}
\label{tab:dich_Bresults}
\begin{tabular}{ccccccc}
\hline
$n$ & $\alpha$ & $\mu$ & $Y$ & $Z_0$ & $Z_1$ & $r$ \\
\hline
$3$  & $6$    & $-\pi/6$    & $4.6\times10^{-14}$ & $0.237$  & $2.1\times10^{2}$ & $6.1\times10^{-14}$ \\
$4$  & $4$    & $-\pi/4$    & $3.1\times10^{-14}$ & $0.158$  & $1.1\times10^{2}$ & $3.8\times10^{-14}$ \\
$5$  & $4$    & $-3\pi/10$  & $3.1\times10^{-14}$ & $0.179$  & $7.0\times10^{1}$ & $3.9\times10^{-14}$ \\
$6$  & $13/5$ & $-\pi/3$    & $2.2\times10^{-14}$ & $0.077$  & $8.6\times10^{1}$ & $2.5\times10^{-14}$ \\
$7$  & $5/2$  & $-5\pi/14$  & $2.2\times10^{-14}$ & $0.075$  & $6.4\times10^{1}$ & $2.5\times10^{-14}$ \\
$10$ & $5/2$  & $-2\pi/5$   & $1.7\times10^{-14}$ & $0.082$  & $3.9\times10^{1}$ & $2.0\times10^{-14}$ \\
$35$ & $5/2$  & $-33\pi/70$ & $1.5\times10^{-14}$ & $0.085$  & $2.5\times10^{1}$ & $1.7\times10^{-14}$ \\
\hline
\end{tabular}
\end{table}
We performed the computations with Julia 1.12.7, \texttt{RadiiPolynomial.jl}~\cite{bib:radiipolynomial} 0.10.0
and \texttt{IntervalArithmetic.jl}~\cite{bib:intervalarithmetic} 1.0.9. Table~\ref{tab:dich_Bresults} shows
the resulting bounds and radii. In all seven cases, $Y$ is of order $10^{-14}$, $Z_0\le0.237$,
and the admissible radius satisfies $r\le6.1\times10^{-14}$. All codes used in the proofs are
available on GitHub~\cite{bib:codes}.
\akt{Our implementation follows the Julia code examples in~\cite{henot2026cap}, which presents a systematic procedure for computer-assisted proofs with \texttt{RadiiPolynomial.jl}.}

\begin{proof}[Proof of Theorem~\ref{thm:main}]
The validated zero supplies a $2$-periodic profile at each prescribed parameter pair.
The coefficient bounds give $|\ba_0|>0.56$ and $|\ba_1|>0.12$, whereas
$r\le6.1\times10^{-14}$. Thus Remark~\ref{rem:dich_closure} gives $\ta_0, \ta_1\ne0$.
Proposition~\ref{prop:dich_star} now gives a star-shaped periodic solution of type
$(p,m)$, with nonconstant modulus and fundamental period $L=2p$, in each of the seven
cases. It remains to verify the origin avoidance and winding number asserted in
Theorem~\ref{thm:main}. Since the center has real coefficients,
\begin{align*}
\operatorname{Re}\sum_{|k|\le N}\ba_ke^{\im\pi kt}&=\ba_0+
\sum_{k=1}^N(\ba_k+\ba_{-k})\cos(\pi kt)\\
&\ge \ba_0-
\sum_{k=1}^N|\ba_k+\ba_{-k}|.
\end{align*}
Since $\nu\ge1$, we have $\sup_t\bigl|\tw(t)-\sum_{|k|\le N}\ba_ke^{\im\pi kt}\bigr|\le
\norm{\ta-\ba}_\nu\le r${, and thus $|\operatorname{Re}\tw(t)-{\operatorname{Re}\sum_{|k|\le N}\ba_ke^{\im\pi kt}}|\le r$ for all $t\in\R$}. Hence the coefficient ball, including its infinite tail, gives
\begin{equation}\label{eq:positive_real_profile}
\operatorname{Re}\tw(t)\ge
\ba_0-\sum_{k=1}^N|\ba_k+\ba_{-k}|-r>0.10
\qquad(t\in\R).
\end{equation}
The last inequality was checked with interval arithmetic in every case. Hence the
orbit avoids the origin, and Remark~\ref{rem:dich_table_correction} gives signed
winding number $-m$.
\end{proof}

\section*{Concluding remarks}\label{sec:conclusion}

We have proved the existence of star-shaped periodic solutions of \eqref{eq:DDE} for seven
parameter pairs $(\alpha,\theta)$ with $\theta=\pi/n$ (Theorem~\ref{thm:main}). These solutions were
observed numerically by Hara~\cite{P2}, and to the best of our knowledge their existence had not
been proved before.

The proof has three steps. First, numerical observations of the Fourier coefficients led to the
ansatz \eqref{eq:main-anzatz}, which writes a solution as a $2$-periodic profile in a rotating
frame and removes the degeneracy that prevents a direct application of the Fourier spectral
method. Second, restricting the profile to real Fourier coefficients and using the resonance
between the delay and the profile period, we fixed the rotation frequency exactly as
$\mu=\theta-\pi/2$ and obtained the real system \eqref{eq:dich_FR}. Third, we validated a zero of
this system with the Radii Polynomial Theorem, using a preconditioned tail so that all Fourier
modes are included in the bounds.

For seven exact parameter pairs we have proved the existence of star-shaped periodic
solutions in the sense of Definition~\ref{def:star_shaped}. In addition to the defining
rotation relation and nonconstant modulus, we have established their fundamental periods,
origin avoidance and signed winding numbers. The proof constructs profiles in the
symmetric real coefficient subspace; no classification of all star-shaped solutions
is asserted.
Theorem~\ref{thm:main} addresses the existence question motivating
Conjecture~\ref{conj:hara} for the class defined here, at finitely many parameter values.
For irrational $\theta/\pi$, Proposition~\ref{prop:dich_star}(c) describes the closure
of a single orbit, conditional on the existence of a periodic rotating-frame profile.
It does not prove existence at irrational parameters or describe all initial histories.
Its annular radii are the extrema of $|w|$, and have not been identified with the radii
in Hara's dynamical conjecture.

Several problems remain toward the resolution of Conjecture~\ref{conj:hara}. First, the
conjecture concerns the dynamics: the instability of the circular orbit $(CO)$ for
$\alpha>\alpha_0$ and the convergence of solutions to the star-shaped periodic solutions. Proving
these statements requires the stability and the attraction of the solutions obtained here, as well
as the instability threshold $\alpha_0$ for $(CO)$. Second, our result holds at finitely many
parameter values. A parameter-dependent validation along $\alpha$ and $\theta$ would extend it to
continuous families and could address the completeness of the periodic families. Third, for
irrational $\theta/\pi$, the conjecture predicts that each orbit fills an annulus with radii
$\sin\theta$ and $1$. Relating this annulus to the extrema of $|w|$ in
Proposition~\ref{prop:dich_star}(c) and proving the existence of the corresponding solutions
are left for future work.

\section*{Acknowledgements}
The authors are grateful to Emiko Ishiwata, Tetsuya Ishiwata and Yukihiko Nakata for
introducing them to the problem of star-shaped periodic solutions and for sharing the
motivation behind it.
SM was supported by JSPS KAKENHI Grant Number 25K17288.
AT was supported by the Japan Science and Technology Agency (JST) FOREST Program (Grant
Number JPMJFR246A), and by JSPS KAKENHI Grant Numbers 24K00538, 25K00922, and 26K00619.

\section*{{Declarations}}
\akt{Generative AI tools, including Claude Code (Anthropic), were used to assist in the derivation and verification of the rigorous bounds, in the implementation and debugging of the computational code, and in the editing of the manuscript text. The authors independently verified all mathematical arguments and rigorous estimates, and the computational results were validated with interval arithmetic. The code is publicly available~\cite{bib:codes}. The authors assume responsibility for all content.}

% **********************************************************
% Appendices
% **********************************************************

% **********************************************************
% Bibliography
% **********************************************************

\bibliographystyle{plain}
\bibliography{papers}
\end{document}